\documentclass[trsc,sglanonrev]{informs4}
\providecommand{\appendix}{\par\setcounter{section}{0}\setcounter{subsection}{0}\renewcommand{\thesection}{Appendix~\Alph{section}}}

\OneAndAHalfSpacedXII 

\usepackage{natbib}
 \bibpunct[, ]{(}{)}{,}{a}{}{,}%
 \def\bibfont{\small}%
\EquationsNumberedThrough    

\TheoremsNumberedThrough     
\ECRepeatTheorems  %

\usepackage{bigstrut,bigdelim,multirow}
\usepackage{booktabs}
\usepackage{array,tabularx}
\usepackage{enumerate}
\usepackage{url}
\usepackage{verbatim}
\usepackage[ruled,vlined,linesnumbered]{algorithm2e}
\SetKwInput{KwInput}{Input}
\SetKwInput{KwOutput}{Output} 
\usepackage{amsmath}
\usepackage{amssymb}
\usepackage{mathrsfs}
\usepackage{longtable}
\usepackage{bm}
\usepackage{lmodern}
\usepackage{xcolor}

\usepackage{appendix}

\usepackage{graphicx}
\usepackage{subcaption}
\usepackage{float}
\usepackage{caption,rotfloat}

\usepackage{xr-hyper}
\usepackage[colorlinks,linkcolor=blue,anchorcolor=blue,citecolor=blue]{hyperref}
\providecommand{\policy}[1]{\textsf{#1}}
\providecommand{\CVaR}{\operatorname{CVaR}}

\begin{document}

\graphicspath{{fig/}{outputs_manuscript/e2_policy_comparison_rerun/figures/}{outputs_manuscript/e3_value_of_information_rerun/figures/}{outputs_manuscript/computational_section_rescue/figures/}}


\RUNAUTHOR{Guo et al.}

\RUNTITLE{Learning and Dynamic Route Recommendation}

\TITLE{Learning and Dynamic Route Recommendation Under Unknown Traveler Compliance and Road Conditions}

\ARTICLEAUTHORS{%
\AUTHOR{Yanru Guo, Siqian Shen, Ruiwei Jiang, Xi Jessie Yang}
\AFF{Department of Industrial and Operations Engineering,
University of Michigan, Ann Arbor, MI 48109,\\
\EMAIL{\{yanruguo, siqian, ruiwei, xijyang\}@umich.edu}}
} 

\ABSTRACT{%
Route planning during evacuations and other large-scale egress operations is inherently dynamic: the planner's recommendations affect future congestion on the roads, influence whether travelers accept their guidance, and determine what the planner can learn from the travelers' actual route choices. In this paper, we consider stochastic dynamic route recommendation to many travelers over a finite horizon under uncertain compliance and road conditions. In each period, a planner assigns routes that lead different origins' traveler groups to safe destinations, and travelers either follow the recommendation or revert to a preferred route. We model compliance using a logistic model and learn its unknown coefficients from observed route travel. We also filter the state of road conditions from noisy telemetry, and use a performance state to capture how the quality of past route recommendation affects travelers' future compliance. To determine the route assignment policy, we formulate a finite-horizon, belief-state dynamic program, show that the information state is sufficient for optimal Markov control, and develop an approximation that combines compliance beliefs with proximal policy optimization (PPO) using a Dirichlet actor over feasible route shares. We conduct numerical studies on instances of Eastern Massachusetts highway network data under different uncertain road conditions and compliance behavior. We demonstrate the efficiency and effectiveness of the PPO method and the value of dynamic adaptive solutions, against static routes. Our solutions reduce the mean and upper-tail costs of all travelers relative to myopic policy benchmarks, even when these benchmarks know the true compliance coefficients. Moreover, updating beliefs and the performance state dynamically also reduces the total cost under accurate, biased, and diffuse priors.
}%



\maketitle

%


\section{Introduction}\label{sec:intro}

Route recommendation for many travelers who share a road network can redistribute demand and ease congestion. This makes it particularly valuable during evacuations, disruptions, and other large-scale egress operations. However, a recommendation is not a flow assignment that travelers always carry out. Travelers may choose to take more familiar or direct routes, and their combined reaction alters the congestion that later traffic faces. A planner that ignores this feedback from user compliance and trust may repeatedly overload attractive routes or issue diversions that few travelers follow. Effective route recommendation must therefore anticipate both network changes and travelers' responses to the information provided \citep{peeta1995multiple,jahn2005system,paz2009behavior}, which remain imperfectly known during routing operations. 

Existing studies show that traveler compliance depends on factors such as relative travel time, route familiarity, prevailing traffic conditions, prior experience, and the perceived credibility of the guidance source \citep{srinivasan2000modeling,chorus2009traveler,sadri2014how}. Moreover, recommendations create a feedback loop. Their realized performance may influence subsequent route acceptance and traveler compliance, and the recommended routes determine which behavioral responses the planner observes and how informative those responses are. In this paper, we consider the optimization of dynamic route planning for multiple travelers at different origins, while learning their compliance attitudes and uncertain road conditions. We seek to reduce the total travel time and penalty for travelers who have not reached safe destinations by the end of the planning horizon.

The literature has studied several core elements related to this problem, but to the best of our knowledge, no studies have integrated them within a unified framework. \citet{guo2026-cdc} formulate a single-period static route assignment model that accounts for uncertain traveler compliance and uses nonlinear objective functions to minimize total travel time under congestion. Their model issues guidance once and cannot update recommendations for later departures as conditions change. Dynamic traffic assignment, evacuation-routing, and behavioral models capture time-dependent congestion, partial control, and uncertain compliance \citep{merchant1978model,paz2009line,samaranayake2018discrete,lessan2022planning,yun2024navigating}, but they do not learn compliance and road conditions while optimizing recommendations over time. 

In this work, we formulate a finite-horizon, belief-state dynamic model. In each period (i.e., decision epoch), a central planner assigns route shares to traveler groups at different origins to direct them to one of multiple destinations (e.g., safe zones). Travelers accept the assigned destination and follow either the recommended route or a destination-consistent default route. To compute a default route, the planner can solve a weighted shortest path problem that accounts for congestion-free travel time and, where available, preferences for fewer turns or toll-free roads. 
We represent traveler compliance by a logistic model with fixed but unknown coefficients that we learn from observed routes taken by the travelers. We can also infer the evolving latent road state from noisy telemetry and use an exponentially smoothed measure of recommendation performance as an operational proxy for route recommendation quality. The planner then incorporates realized route choices into a dynamic traffic-loading process and minimizes expected arc-entry travel cost together with terminal penalties for unsatisfied demand and residual traffic.

We describe the main contributions of the paper as follows. First, our dynamic model integrates destination-consistent default routes with compliance learning, road state filtering, and recommendation performance feedback. Under mild assumptions, we establish that the resulting belief state is sufficient for optimal Markov control. Second, we develop a computational approximation that represents compliance beliefs with weighted particles and uses a PPO actor with Dirichlet distributions to generate feasible route-share allocations. Third, we evaluate dynamic control and approximation accuracy on a small instance for which exhaustive enumeration gives the minimum cost over a finite set of route-share allocations, and assess mean cost, upper-tail cost, and the joint value of online adaptation on moderate, real-world sized instances based on the Eastern Massachusetts highway network data. 

The remainder of the paper is organized as follows. In Section~\ref{sec:review}, we review the relevant literature and position our contributions. In Section~\ref{sec:problem-description}, we describe our problem setting and introduce the notation. In Section~\ref{sec:dynamic-route-model}, we formulate the stochastic dynamic program and Bellman recursion for optimizing dynamic route recommendation. In Section~\ref{sec:ppo-approximation}, we develop the belief representation and approximation method. We conduct numerical studies and report results in Section~\ref{sec:numerical}. In Section~\ref{sec:conclusion}, we conclude the work and discuss future research directions.

\section{Literature Review}\label{sec:review}

Our work draws on four research streams: dynamic route recommendation and evacuation routing, traveler responses to guidance, Bayesian control under partial observability, and approximate dynamic programming and reinforcement learning in transportation.

Dynamic traffic assignment (DTA; \citealp{merchant1978model}) represents time-varying demand, traffic propagation, and congestion, and it distinguishes user-equilibrium assignments, which reflect travelers' individual route choices, from system-optimal assignments that minimize network-wide costs \citep{zheng2015system}. \citet{peeta1995multiple} examine rolling-horizon methods that update assignments as network conditions evolve. Other models let travelers retain some control over their routes. Constrained system-optimal models limit excessive detours \citep{jahn2005system}, behaviorally consistent guidance anticipates traveler responses \citep{paz2009behavior}, and \citet{samaranayake2018discrete} study system-optimal DTA when the planner controls only a fraction of total flow. \citet{paz2009line} recalibrate behavioral parameters online. Evacuation routing adds the need for rapid network clearance across safe destinations with imbalanced roadway capacity. Existing studies jointly optimize destination and route choices \citep{han2006global}, incorporate departure-time decisions \citep{chiu2007modeling}, and reroute traffic in real time under road disruptions \citep{darvishan2021dynamic}. Recent work incorporates imperfect compliance into dynamic \citep{lessan2022planning,yu2026multi} or static route planning \citep{guo2026-cdc}. However, these studies typically represent compliance and disruptions with fixed parameters or predefined scenarios rather than with beliefs that the planner updates within a stochastic dynamic optimization framework, as we do in this work.

Travelers balance route recommendations against a preference for familiar routes, and information quality, network conditions, prior experience, and switching costs influence this balance \citep{srinivasan2000modeling, jia2020context}. Experimental studies associate compliance with the characteristics of a recommendation, congestion or incidents, and prior experience \citep{chen2003driver}. Familiarity matters especially during evacuations, when travelers may prefer habitual routes to recommended alternatives \citep{sadri2014how}. Traveler responses also evolve through experience. Travelers combine real-time information with personal observations \citep{ben2008combined}, account for the reliability of the information source \citep{chorus2009traveler,ben2013impact}, and may reduce their reliance after poor recommendations \citep{fox1998effects}. These findings motivate our state-dependent compliance model and its recommendation performance feedback.

Partially observable Markov decision processes (POMDPs) distinguish among latent states, observations, and control actions \citep{smallwood1973optimal}, and Bayesian filtering produces a belief state that is sufficient for dynamic programming \citep{monahan1982state}. When parameters are unknown, Bayes-adaptive POMDPs augment the state with a posterior distribution over them \citep{ross2011bayesian}. Control actions then influence both immediate performance and the information available for future decisions, which creates a dual-control effect \citep{klenske2016dual}. In transportation, partially observable control has mainly addressed traveler-level adaptation, such as routing with stochastic travel times \citep{gao2006optimal, gao2012real, huang2018trajectory}. Approximate dynamic programming \citep{powell2010approximate} and reinforcement learning \citep{sutton2000policy} support mobility-on-demand operations \citep{kullman2022dynamic} and adaptive routing \citep{mao2018reinforcement, jia2020context}, and \citet{yun2024navigating} account for traveler deviations from system-oriented route recommendations. These studies generally represent behavioral uncertainty through prespecified response rules or scenarios rather than through a posterior distribution that the planner updates, as we do in this work.

\section{Problem Description}
\label{sec:problem-description}
Consider an emergency manager who must move thousands of residents out of a threatened area within a few hours. The manager can tell each resident which route to take but cannot make anyone take it. Residents who doubt the guidance, or who simply prefer a road they know, fall back on a familiar route. Because familiar routes often converge on the same few corridors, this fallback behavior can recreate the congestion that the guidance sought to avoid. At the same time, an incident on a key road can change which routes are attractive, and the manager typically learns about it only through noisy reports. In this section, we formalize the decision problem that such a planner faces.

We consider route recommendation for evacuations and other large-scale egress operations. A central planner periodically issues route recommendations to travelers from a finite set of traveler groups. Each group departs from an origin node, and multiple groups may share the same origin. The planner allocates travelers from each group to routes that link their origin with one of their acceptable destinations. Because the planner cannot control travelers after departure, we feed their realized route choices into a dynamic traffic-loading process that propagates congestion over time. The same framework applies to campus and industrial emergency response, large-event egress, and other settings in which multiple exits, shelters, assembly areas, or mobility hubs serve as acceptable destinations.

Let $G=(\mathcal V,\mathcal E)$ denote a directed transportation network, where $\mathcal V$ is the set of nodes and $\mathcal E\subseteq\mathcal V\times\mathcal V$ is the set of directed arcs. Let $\mathcal O$ be a finite index set of traveler groups, and let $v:\mathcal O\rightarrow\mathcal V$ be a mapping, where $v(o)$ denotes the origin associated with group $o$. Distinct groups $o\neq o'$ may share an origin, that is, $v(o) = v(o')$. Separating groups from physical locations lets the planner treat co-located populations differently, for example, residents and visitors who know the local roads to different degrees. 
Let $\mathcal D\subseteq\mathcal V$ contain all destinations. For each group $o\in\mathcal O$, the nonempty set $\mathcal D_o\subseteq\mathcal D$ contains the destinations considered safe to travelers in that group. A route for group $o$ to destination $d\in\mathcal D_o$ is an ordered sequence of consecutive directed arcs that begins at origin $v(o)$ and terminates at $d$.

The planning horizon consists of decision epochs $t=0,1,\ldots,T-1$ with period length $\Delta t>0$, followed by the terminal time $T$. Travel time and $\Delta t$ share the same unit, and we express arc capacities as numbers of travelers per decision period. In the baseline model, every destination in $\mathcal D_o$ has sufficient receiving capacity for the modeled population, and travelers accept any destination in $\mathcal D_o$ that the planner assigns.
(Binding destination capacities would require an extended state and an assignment mechanism that respects the remaining receiving capacities.)
When $|\mathcal D_o|>1$, the planner jointly chooses the destination and the route for travelers in group $o$. When $\mathcal D_o=\{d_o\}$, the destination is fixed, and the planner selects among alternative routes to $d_o$. Fixed-destination route recommendation is therefore a special case of the destination-and-route recommendation model that we study. 

\subsection{Recommendations and Destination-Consistent Default Routes}
\label{subsec:problem-setting}

At each decision epoch, the planner splits each group's departing travelers among the routes to the group's acceptable destinations.
We assume that travelers always accept their assigned destinations but may decline the recommended route. A noncompliant traveler follows a predetermined default route to the assigned destination. In the baseline model, this default route is the nominal free-flow shortest path, which the planner can compute \emph{a priori}.
The default route serves as a proxy for the familiar or direct route to which travelers revert when they reject guidance \citep{srinivasan2000modeling,sadri2014how}. As a result, the planner never faces an arbitrary deviation. Each traveler makes a binary choice between the recommendation and a predictable fallback, and this structure makes compliance both learnable and plannable.
(Online Supplement Section~\ref{app:multiple-default-routes} extends the model to let noncompliant travelers select among multiple destination-consistent alternatives.)

\begin{remark}
\label{remark:1}
    In the baseline model, each assigned traveler has two possible route outcomes: the recommended route and the corresponding default route. When the two routes coincide, the planner observes the realized path but cannot tell whether the traveler intended to follow the guidance. This distinction matters in our later analysis. The same realized route counts drive both traffic loading and statistical learning, but only recommendations that differ from the default routes reveal information about compliance. From a managerial perspective, a recommendation that agrees with travelers' habits is easy to issue but teaches the planner nothing. A recommendation that departs from habit carries a learning value in addition to its effect on congestion.
\end{remark}

\subsection{Information and Decision Timeline}
\label{subsec:problem-learning-performance}

The planner works in a repeated cycle of deciding, observing, and learning. Let $S_t$ be the state variable that collects the information available immediately before the planner selects the period-$t$ route recommendations. Section~\ref{subsec:dp-state} defines the components of $S_t$. At the start of period $t$, the planner observes $S_t$, which includes waiting demand, the inherited traffic state, beliefs about compliance and road conditions, and the recommendation performance state. The planner then selects route shares for each traveler group. Next, the demand process determines the actual departure counts, the planner assigns individual travelers to routes according to the shares, and travelers follow either their recommended routes or their default routes as traffic evolves. Finally, new demand joins the origin queues, and the planner receives route agreement observations and road condition telemetry, updates its beliefs and the recommendation performance state, and enters period $t+1$ with state $S_{t+1}$.

\begin{figure}[htbp!]
  \centering
  \includegraphics[width=0.75\linewidth]{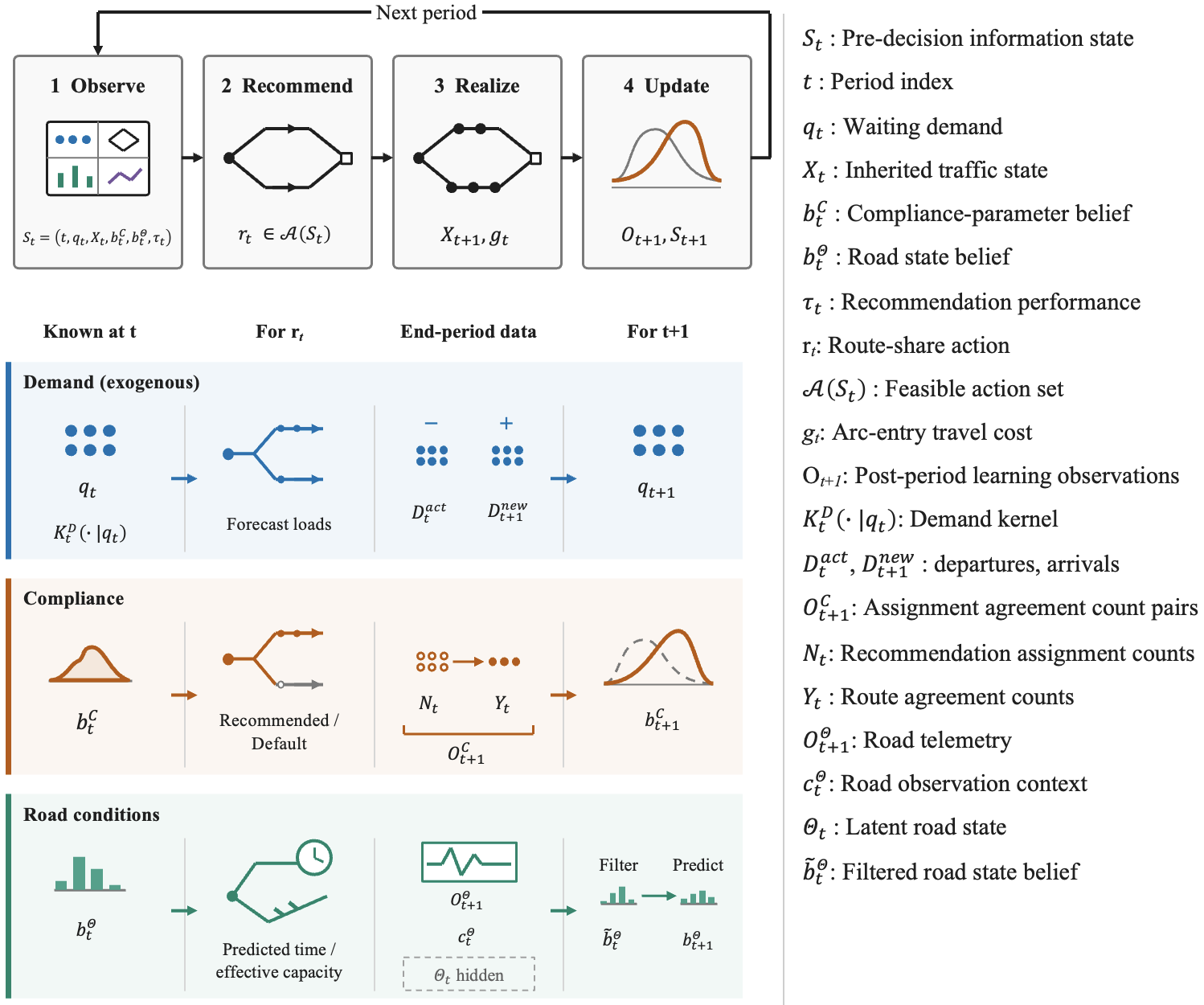}
  \caption{Within-period information and decision timeline for dynamic route recommendation. The planner observes $S_t$, selects route shares $\mathbf r_t$, and updates $S_{t+1}$ after traffic realization and information updates. The three rows trace exogenous demand, uncertain compliance, and latent road conditions.}
  \label{fig:route-guidance-timeline}
\end{figure}

Three interdependent sources of uncertainty shape this cycle: (i) actual departures may differ from forecasts, (ii) travelers may not follow recommended routes, and (iii) the road state may be unobserved and may change because of emergencies and disruptions. The planner learns about (ii) and (iii) online. In addition, recent recommendations may influence subsequent compliance through an endogenous performance state. This feedback turns the credibility of the guidance into a resource that the planner builds or spends with every recommendation. A recommendation that performs reliably may raise compliance in later periods, whereas a diversion that turns out slower than the travelers' default route may erode the planner's ability to steer later departures. 
%
Figure~\ref{fig:route-guidance-timeline} depicts the timeline and shows how demand, compliance, and road conditions inform route recommendations and the updates that lead to $S_{t+1}$. Sections~\ref{sec:dynamic-route-model}--\ref{sec:ppo-approximation} define the detailed notation. 

\subsection{Objective}
\label{subsec:problem-objective}

Let $\mu=(\mu_0,\ldots,\mu_{T-1})$ denote an admissible dynamic route recommendation policy, whose period-$t$ decision rule uses the information available before that decision to select a feasible route-share action $\mathbf r_t$. We seek a policy that solves
\begin{equation}
\label{eq:problem-objective}
\min_{\mu}
\mathbb E^{\mu}\left[
\sum_{t=0}^{T-1} g_t(S_t,\mathbf r_t,\Xi_{t+1})
+G_T(S_T)
\mid S_0\right], 
\end{equation}
where $\mathbb E^{\mu}$ averages over trajectories generated by policy $\mu$ conditional on the initial state $S_0$, $\Xi_{t+1}$ collects the within-period random variables, $g_t$ is the arc-entry travel cost in period $t$, and $G_T$ penalizes travelers still waiting at their origins or still en route at time $T$. The objective reflects two concerns of an emergency manager. The stage cost measures the travel time that the population experiences, and the terminal penalty measures how many people the operation fails to bring to safety within the horizon. Because the terminal penalty depends on the random terminal state, the policy affects the objective through both its current actions and the resulting evolution of future states. A planner who minimizes only the current travel cost can therefore leave many travelers en route at time $T$. Section~\ref{sec:numerical} shows that this trade-off drives most of the performance differences among the policies we compare. 

\section{Belief-State Dynamic Route Recommendation Model}
\label{sec:dynamic-route-model}

We now turn the description in Section~\ref{sec:problem-description} into a dynamic program and follow the order of the timeline in Figure~\ref{fig:route-guidance-timeline}. Section~\ref{subsec:dp-state} defines the state, routes, and feasible actions. Sections~\ref{subsec:dp-demand-counts}--\ref{subsec:dp-traffic-cost} trace the events within a period: departures, compliance responses and their Bayesian update, road state filtering, and traffic propagation, performance feedback, and costs. Section~\ref{subsec:dp-transition} assembles the controlled transition law, shows that the belief state suffices for optimal control, and states the Bellman recursion.

\subsection{State, Routes, and Feasible Actions}
\label{subsec:dp-state}

Let $\mathscr F_t$ be the $\sigma$-field generated by the initial information, previous actions, and observations received through the end of period $t-1$. Rather than carry this growing history, the planner summarizes it in the state 
\begin{equation}
\label{eq:dp-state}
    S_t
    =
    \left(
        t,
        \mathbf q_t,
        \mathbf X_t,
        b_t^C,
        b_t^\Theta,
        \boldsymbol\tau_t
    \right).
\end{equation}
Here the queue vector $\mathbf q_t=(q_{o,t})_{o\in\mathcal O}$ counts travelers waiting in each group. Groups at the same physical location keep separate queues because they may have different acceptable destination sets.
The traffic state $\mathbf X_t$ tracks travelers who have departed but have not reached their assigned safe destinations. In our traffic-loading model, a packet is a group of travelers with the same realized route, position along that route, and remaining traversal duration. Let $C_t$ be the number of packets at the beginning of period $t$. The traffic state is the collection of packet records 
\begin{equation}
\label{eq:traffic-packet-state}
    \mathbf X_t = \bigl(n_{c,t},p_{c,t},k_{c,t},\ell_{c,t}\bigr)_{c=1}^{C_t}.
\end{equation}
For packet $c$, the size $n_{c,t}\in\mathbb Z_{\geq1}$ counts its travelers, and the route $p_{c,t}$ contains $|p_{c,t}|$ arcs. The index $k_{c,t}\in\{1,\ldots,|p_{c,t}|\}$ identifies the packet's most recently entered arc, and $\ell_{c,t}\in\mathbb Z_{\geq0}$ is the number of periods remaining on that arc. A packet with $\ell_{c,t}=0$ and $k_{c,t}<|p_{c,t}|$ enters its next arc during period $t$. Each packet retains its realized route while it remains in the network. At each period, we remove packets that have completed their last arc and reindex the remaining records with $c$ when forming $\mathbf X_t$.
The planner observes all packet records at each period. The collection is empty when $C_t=0$. Otherwise, each record contains a route identifier and three integers. Let $e_{p,k}$ be the $k^\text{th}$ arc of route $p$. The packets that are ready to move determine the inherited arc-entry count 
\begin{equation}
\label{eq:inherited-entry-flow}
    h_{ij,t}(\mathbf X_t)
    =
    \sum_{\substack{1\leq c\leq C_t\\
                    \ell_{c,t}=0,\ k_{c,t}<|p_{c,t}|}}
    n_{c,t}\mathbf 1\{e_{p_{c,t},k_{c,t}+1}=(i,j)\},
    \quad (i,j)\in\mathcal E.
\end{equation}
This quantity counts entries by previously departed travelers, not current arc occupancy. It tells the planner how much of each arc's period-$t$ capacity earlier decisions have already committed, and therefore how much room remains for new departures.

A fixed but unknown coefficient vector $\boldsymbol\beta^C\in\mathbb R^5$ governs compliance and applies across traveler groups, routes, destinations, and periods. Its entries are the (i) intercept and coefficients on (ii) scaled route-time difference, (iii) familiarity, (iv) recommendation performance, and (v) congestion. Because the groups share this vector, every informative recommendation, wherever the planner issues it, sharpens the planner's view of all travelers. The compliance belief is the joint posterior measure
\begin{equation}
\label{eq:compliance-belief}
    b_t^C(d\boldsymbol\beta^C)
    =
    \Pr\left(
        \boldsymbol\beta^C\in d\boldsymbol\beta^C
        \mid \mathscr F_t
    \right).
\end{equation}

Let $\Theta_t\in\mathcal S^\Theta$ denote the road state, where $\mathcal S^\Theta$ is a known finite set of possible road conditions. Its belief is the probability vector $b_t^\Theta=(b_t^\Theta(\vartheta))_{\vartheta\in\mathcal S^\Theta}$, with entries
\begin{equation}
\label{eq:road-belief}
    b_t^\Theta(\vartheta)
    =
    \Pr(\Theta_t=\vartheta\mid\mathscr F_t),
    \quad \vartheta\in\mathcal S^\Theta.
\end{equation}
The entries lie in $[0,1]$ and sum to one, and thus the belief has $|\mathcal S^\Theta|-1$ independent coordinates. Finally, the vector $\boldsymbol\tau_t=(\tau_{o,t})_{o\in\mathcal O}\in[0,1]^{|\mathcal O|}$ records recommendation performance, an operational proxy for the quality of recent guidance as travelers experience it. Section~\ref{subsec:dp-traffic-cost} defines its snapshot-based performance criterion and exponential smoothing update.

The initial physical state $(\mathbf q_0,\mathbf X_0)$, beliefs $b_0^C,b_0^\Theta$, and performance vector $\boldsymbol\tau_0\in[0,1]^{|\mathcal O|}$ are model inputs. We assume that $\boldsymbol\beta^C$ and $\Theta_0$ are conditionally independent given $\mathscr F_0$. That is, for every measurable $B\subseteq\mathbb R^5$ and $\vartheta\in\mathcal S^\Theta$, we have 
\begin{equation}
\label{eq:initial-belief-factorization}
    \Pr\!\left(
        \boldsymbol\beta^C\in B,\Theta_0=\vartheta
        \mid\mathscr F_0
    \right)
    =
    b_0^C(B)b_0^\Theta(\vartheta).
\end{equation}
We define the set of active traveler groups as $\mathcal O_t^+(S_t) = \{o\in\mathcal O:q_{o,t}>0\}$. Because actual departures satisfy $D^{\mathrm{act}}_{o,t}\leq q_{o,t}$, only groups in $\mathcal O_t^+(S_t)$ can generate departures and require action coordinates during period $t$. For each group $o$ and acceptable destination $d\in\mathcal D_o$, let $\mathcal P_{od,t}(S_t)$ be the finite set of routes known to be available from physical origin $v(o)$ to $d$. The group route set is
$\mathcal P_{o,t}(S_t)=\bigcup_{d\in\mathcal D_o}\mathcal P_{od,t}(S_t)$, which is nonempty for every active group. Each route $p$ terminates at a unique destination $d(p)\in\mathcal D_o$.
Using the known nominal free-flow arc time $\bar t^0_{ij,t}$, we define the nominal route travel time $\operatorname{TT}^0_{p,t} =\sum_{(i,j)\in p}\bar t^0_{ij,t}$.
In the experiments, the nominal travel time equals the normal-state value of the free-flow functions in~\eqref{eq:dp-travel-time}. The nominal travel time and route catalog are time-invariant, and we suppress $t$ when reporting the free-flow travel time in Online Supplement Section~\ref{app:ema-implementation}. 
For each nonempty $\mathcal P_{od,t}(S_t)$, we define the destination-specific default route by
\begin{equation}
\label{eq:destination-default-route}
    p^0_{od,t}(S_t)
    \in
    \arg\min_{q\in\mathcal P_{od,t}(S_t)}
    \operatorname{TT}^0_{q,t},
\end{equation}
where we break ties by a fixed deterministic rule. For a recommended route $p$, we set $p^0_{op,t}:=p^0_{o,d(p),t}(S_t)$ as its destination-consistent default route. This construction ensures that a traveler who rejects a recommendation still heads to the destination that the planner assigned, and thus noncompliance changes the path but not the destination.

The planner's period-$t$ decision is the route-share action
\begin{equation}
\label{eq:dp-action}
    \mathbf r_t
    =
    \left\{
        r_{op,t}:
        o\in\mathcal O_t^+(S_t),
        \ p\in\mathcal P_{o,t}(S_t)
    \right\}.
\end{equation}
We write $\mathbf r_{o,t}=(r_{op,t})_{p\in\mathcal P_{o,t}(S_t)}$ for group $o$'s share vector. Given state $S_t$, the feasible action set is
\begin{equation}
\label{eq:dp-action-set}
    \mathcal A(S_t)
    =
    \left\{
        \mathbf r_t:
        \mathbf r_{o,t}\in
        \mathbb R_{\geq0}^{|\mathcal P_{o,t}(S_t)|},
        \sum_{p\in\mathcal P_{o,t}(S_t)}r_{op,t}=1,
        \forall o\in\mathcal O_t^+(S_t)
    \right\}.
\end{equation}
When no groups are active, $\mathcal A(S_t)$ contains a single action, which issues no recommendations. Route shares, rather than individual assignments, are the natural decision for a planner who must commit before knowing how many travelers will actually depart. The same share vector applies to any realized departure count.

\subsection{Exogenous Demand and Recommendation Counts}
\label{subsec:dp-demand-counts}

Departures are the first random event in each period, and the planner cannot influence them. Let $\mathbf D_t^{\mathrm{act}} =(D_{o,t}^{\mathrm{act}})_{o\in\mathcal O}$ and $\mathbf D_{t+1}^{\mathrm{new}} =(D_{o,t+1}^{\mathrm{new}})_{o\in\mathcal O}$ be the vectors of departures from origin queues during period $t$ and new travelers who join those queues after period $t$'s traffic loading, respectively. New travelers become eligible to depart in period $t+1$. The planner observes both $\mathbf D_t^{\mathrm{act}}$ and $\mathbf D_{t+1}^{\mathrm{new}}$ at the end of period $t$. In the baseline model, the conditional distribution of $(\mathbf D_t^{\mathrm{act}},\mathbf D_{t+1}^{\mathrm{new}})$ given the observable history $\mathscr F_t$ depends only on $t$ and $\mathbf q_t$, through the known time-indexed kernel $K_t^D(\cdot\mid\mathbf q_t)$. Its support consists of nonnegative integer vectors $(\mathbf d^{\mathrm{act}},\mathbf d^{\mathrm{new}})$ satisfying $0\leq d^{\mathrm{act}}_o\leq q_{o,t}$ and $d^{\mathrm{new}}_o\geq0$ for every $o\in\mathcal O$. We assume that further conditioning on the selected action $\mathbf r_t$, the compliance coefficients $\boldsymbol\beta^C$, and the current road state $\Theta_t$ leaves this demand law unchanged. Thus, recommendations shape where travelers go but not when they leave.
After the demand process determines the departure counts, the planner randomizes individual assignments according to the selected shares, which yields the recommendation counts
\begin{equation}
\label{eq:recommendation-multinomial}
    \left(N_{op,t}\right)_{p\in\mathcal P_{o,t}(S_t)}
    \mid S_t, \mathbf D_t^{\mathrm{act}},\mathbf r_{t}
    \sim
    \operatorname{Multinomial}\!\left(
        D^{\mathrm{act}}_{o,t},\mathbf r_{o,t}
    \right),
    \quad o\in\mathcal O_t^+(S_t).
\end{equation}
Conditional on $S_t$, $\mathbf D_t^{\mathrm{act}}$, and $\mathbf r_t$, the assignment vectors are independent across traveler groups. Further conditioning on $\mathscr F_t$, $\boldsymbol\beta^C$, $\Theta_t$, and $\mathbf D_{t+1}^{\mathrm{new}}$ leaves their conditional laws and independence unchanged. For every traveler group $o$, $N_{op,t}\in\mathbb Z_{\geq0}$ and $\sum_{p\in\mathcal P_{o,t}(S_t)}N_{op,t} =D^{\mathrm{act}}_{o,t}$, with conditional means $D^{\mathrm{act}}_{o,t}r_{op,t}$. Inactive groups have zero departures and generate no assignment counts. Let $\mathbf N_t$ collect the assignment counts of all active groups. The same realized counts drive both traffic loading and statistical updating. This coupling is the source of a dual-control effect. The shares that the planner chooses to manage congestion also determine how much it learns about travelers.

\subsection{Compliance Response and Bayesian Updating}
\label{subsec:dp-compliance-learning}

The planner cannot observe a traveler's intent, but it can predict how likely a traveler is to follow a given recommendation. For each $o\in\mathcal O_t^+(S_t)$ and $p\in\mathcal P_{o,t}(S_t)$, the compliance feature vector is
\begin{equation}
\label{eq:compliance-feature-vector}
    \boldsymbol\varphi_{op,t}
    =
    \left(
        1,
        \Delta\mathrm{TT}_{op,t},
        \mathrm{Fam}_{op},
        \tau_{o,t},
        \mathrm{Cong}_{op,t}
    \right)^\top
    \in\mathbb R^5.
\end{equation}
We set the first component to $1$ so that its coefficient serves as the intercept, contributing an additive constant to the linear predictor $(\boldsymbol\beta^C)^\top\boldsymbol\varphi_{op,t}$.
The remaining four components capture the recommended route's predicted travel time relative to its default, the group's familiarity with the route, the performance of previous recommendations, and current route congestion.

Let $\widehat{\operatorname{TT}}^{\mathrm{pred}}_{p,t}$ denote the pre-decision prediction of travel time on route $p$, which the planner computes from the current state $S_t$, including the road state belief $b_t^\Theta$. With a scale $s_{\mathrm{TT}}>0$ in the same time units, we define 
\begin{equation}
\label{eq:compliance-travel-time-difference}
    \Delta\mathrm{TT}_{op,t}
    =
    \frac{
        \widehat{\operatorname{TT}}^{\mathrm{pred}}_{p,t}
        -
        \widehat{\operatorname{TT}}^{\mathrm{pred}}_{p^0_{op,t},t}
    }{s_{\mathrm{TT}}}, 
\end{equation}
such that positive $\Delta\mathrm{TT}_{op,t}$ indicates a longer predicted travel time for the recommended route than for its destination-consistent default.

In \eqref{eq:compliance-feature-vector}, the familiarity feature $\mathrm{Fam}_{op}\in[0,1]$ is a fixed, known score of group $o$'s familiarity with route $p$, and larger values indicate greater familiarity. The performance feature $\tau_{o,t}\in[0,1]$ is the state component in~\eqref{eq:dp-state}, which reflects feedback received before the current decision. (We give its recursive update in~\eqref{eq:performance-ema}.) Finally, let $\widehat c_{ij,t}>0$ be the belief-based effective capacity for period $t$, expressed in the same units as the inherited entry count $h_{ij,t}(\mathbf X_t)$. We define
\begin{equation}
\label{eq:compliance-congestion-feature}
    \mathrm{Cong}_{op,t}
    =
    \min\!\left\{
        C_{\mathrm{Cong}},
        \frac{1}{|p|}
        \sum_{(i,j)\in p}
        \frac{h_{ij,t}(\mathbf X_t)}{\widehat c_{ij,t}}
    \right\},
    \quad C_{\mathrm{Cong}}>0.
\end{equation}
This feature is the average inherited entry-to-capacity ratio along $p$, capped at $C_{\mathrm{Cong}}$, where~\eqref{eq:inherited-entry-flow} defines the inherited entry count $h_{ij,t}$. The route-time and effective-capacity maps are known measurable functions of $S_t$. All five features are therefore $\mathscr F_t$-measurable and remain fixed when we evaluate the compliance likelihood. Together, they represent the roles of relative route performance, familiarity, prior experience, and perceived credibility in travelers' responses to guidance \citep[see][]{chen2003driver,ben2008combined,chorus2009traveler,srinivasan2000modeling,sadri2014how}. The performance and congestion features depend on the planner's earlier recommendations, and thus today's recommendations shape tomorrow's compliance.

A traveler assigned to $p$ follows either $p$ or its destination-consistent default $p^0_{op,t}$ and keeps destination $d(p)$ in either case. Let $Y_{op,t}$ count travelers assigned to $p$ whose selected route agrees with the recommendation, and let $\mathbf Y_t$ collect these counts over the same cells as $\mathbf N_t$. The informative recommendation cells are
$\mathcal P^{\mathrm{inf}}_{o,t}(S_t)=\{p\in\mathcal P_{o,t}(S_t):p\neq p^0_{op,t}\}$.
For each informative cell, the probability of agreement conditional on $\boldsymbol\beta^C$ is
\begin{equation}
\label{eq:logistic-compliance}
    \pi^C_{op,t}
    =
    \sigma\!\left(
        (\boldsymbol\beta^C)^\top\boldsymbol\varphi_{op,t}
    \right),
    \quad
    \sigma(x)=\frac{1}{1+e^{-x}},
\end{equation}
and the agreement count follows
\begin{equation}
\label{eq:compliance-binomial}
    Y_{op,t}
    \mid
    \mathbf N_t,\boldsymbol\beta^C,S_t,\mathbf r_t
    \sim
    \operatorname{Binomial}\!\left(
        N_{op,t},\pi^C_{op,t}
    \right),
    \quad
    p\in\mathcal P^{\mathrm{inf}}_{o,t}(S_t).
\end{equation}
If the recommendation and the default coincide, both responses produce the same route, and 
\begin{equation}
\label{eq:coincident-route-agreement}
    Y_{op,t}=N_{op,t}
    \quad\text{almost surely},
    \quad
    \forall p\in\mathcal P_{o,t}(S_t)
    \setminus\mathcal P^{\mathrm{inf}}_{o,t}(S_t).
\end{equation}
These cells contribute to traffic loading but reveal nothing about $\boldsymbol\beta^C$, as Remark~\ref{remark:1} anticipates.

Conditional on $\boldsymbol\beta^C$, $S_t$, $\mathbf r_t$, and $\mathbf N_t$, the agreement counts are mutually independent, and further conditioning on $\mathscr F_t$, $\Theta_t$, and $\mathbf D_{t+1}^{\mathrm{new}}$ leaves this response law and conditional independence unchanged. (Online Supplement Section~\ref{app:beta-binomial} provides an extension with over-dispersion.) 

In the baseline model, travelers commit to complete routes at departure and report their selected routes, matched to their recommendations, by the end of that period. The planner thus receives
$   O^C_{t+1}
    =
    \left\{
        (N_{op,t},Y_{op,t}):
        o\in\mathcal O_t^+(S_t),
        \ p\in\mathcal P_{o,t}(S_t)
    \right\}.$ 
Conditional on the assignment counts, the compliance likelihood is
\begin{equation}
\label{eq:compliance-likelihood}
    L_t^C(\boldsymbol\beta^C)
    =
    \prod_{o\in\mathcal O_t^+(S_t)}
    \prod_{p\in\mathcal P^{\mathrm{inf}}_{o,t}(S_t)}
    \binom{N_{op,t}}{Y_{op,t}}
    (\pi^C_{op,t})^{Y_{op,t}}
    (1-\pi^C_{op,t})^{N_{op,t}-Y_{op,t}},
\end{equation}
where we evaluate each $\pi^C_{op,t}$ at the candidate $\boldsymbol\beta^C$ using~\eqref{eq:logistic-compliance}. Once the pre-decision information and action are fixed, the assignment probabilities contain no unknown compliance parameters, and their factors cancel in Bayes' rule. Under the baseline independence assumptions, the exact update is
\begin{equation}
\label{eq:bayes-compliance-update}
    b_{t+1}^C(d\boldsymbol\beta^C)
    =
    \frac{
        L_t^C(\boldsymbol\beta^C)
        b_t^C(d\boldsymbol\beta^C)
    }{
        \displaystyle
        \int_{\mathbb R^5}
        L_t^C(\widetilde{\boldsymbol\beta}^C)
        b_t^C(d\widetilde{\boldsymbol\beta}^C)
    }.
\end{equation}
If no informative cell has a positive assignment count, then $L_t^C\equiv1$ and $b_{t+1}^C=b_t^C$. In other words, the planner learns about compliance only when it asks some travelers to deviate from their default routes. Online Supplement Section~\ref{app:compliance-posterior} derives this posterior update.  

\subsection{Hidden Road State and Filtering}
\label{subsec:dp-road-belief}

Compliance is not the only unknown. The unobserved road state $\Theta_t$ determines the period-$t$ travel-time and capacity functions, and an incident on a single critical arc can change which routes are attractive. We specify the road transition probabilities as model inputs, and the planner infers the current road state from period-end observations. The road process is exogenous. Conditional on $\Theta_t$, its next state is independent of the compliance coefficients and of all information received through the end of period $t$, that is, 
\begin{equation}
\label{eq:road-transition-kernel}
    \Pr\!\left(\Theta_{t+1}=\vartheta'
        \mid\Theta_t=\vartheta,\boldsymbol\beta^C,\mathscr F_{t+1}\right)
    =\Pr(\Theta_{t+1}=\vartheta'\mid\Theta_t=\vartheta).
\end{equation}
After traffic evolves, the planner receives a joint road observation $O^\Theta_{t+1}$ that contains all retained measurements informative about $\Theta_t$, or a sufficient statistic for them. These measurements include the road state information carried by observed packet transitions, entry durations, and stage costs. The law of $O^\Theta_{t+1}$ depends on the action and realized flows through the observed context
\begin{equation}
\label{eq:road-filter-context}
    \mathbf c_t^\Theta
    =
    \left(
        S_t,
        \mathbf r_t,
        \mathbf D_t^{\mathrm{act}},
        \mathbf N_t,
        \mathbf Y_t
    \right).
\end{equation}
For a measurable set $A$ in the telemetry space, let
\begin{equation}
\label{eq:road-observation-kernel}
    Q_t^\Theta(A\mid\vartheta,\mathbf c_t^\Theta)
    =
    \Pr\!\left(
        O^\Theta_{t+1}\in A
        \mid
        \Theta_t=\vartheta,
        \mathbf c_t^\Theta
    \right).
\end{equation}
We assume that further conditioning on $\boldsymbol\beta^C$, $\mathbf D_{t+1}^{\mathrm{new}}$, and the pre-decision history $\mathscr F_t$ leaves this observation law unchanged. With a common dominating measure
$\Lambda_t(\cdot\mid\mathbf c_t^\Theta)$, we define the likelihood
\begin{equation}
\label{eq:road-observation-likelihood}
    L_t^\Theta(\xi^\Theta\mid\vartheta,\mathbf c_t^\Theta)
    =
    \frac{
        dQ_t^\Theta(\cdot\mid\vartheta,\mathbf c_t^\Theta)
    }{
        d\Lambda_t(\cdot\mid\mathbf c_t^\Theta)
    }(\xi^\Theta).
\end{equation}

Under the baseline assumptions, the counts in $\mathbf c_t^\Theta$ carry no additional information about $\Theta_t$ given $\mathscr F_t$ and $\mathbf r_t$. Let $\widetilde b_t^\Theta$ be the posterior for $\Theta_t$ after the period-$t$ observations. The exact current-state filtering and next-state prediction steps are
\begin{align}
\label{eq:road-filter-update}
    \widetilde b_t^\Theta(\vartheta)
    & =
    \frac{
        L_t^\Theta(
            O^\Theta_{t+1}\mid
            \vartheta,\mathbf c_t^\Theta
        )
        b_t^\Theta(\vartheta)
    }{
        \displaystyle
        \sum_{u\in\mathcal S^\Theta}
        L_t^\Theta(
            O^\Theta_{t+1}\mid
            u,\mathbf c_t^\Theta
        )
        b_t^\Theta(u)
    }, \\
\label{eq:road-filter-predict}
    b_{t+1}^\Theta(\vartheta')
    & =
    \sum_{\vartheta\in\mathcal S^\Theta}
    \Pr(\Theta_{t+1}=\vartheta'\mid\Theta_t=\vartheta)
    \widetilde b_t^\Theta(\vartheta), 
\end{align}
where \eqref{eq:road-filter-predict} uses the transition probabilities in~\eqref{eq:road-transition-kernel}.
For observations with positive predictive likelihood, \eqref{eq:road-filter-update} converts the prior $b_t^\Theta$ into $\widetilde b_t^\Theta$, the posterior distribution of the completed period's road state $\Theta_t$. Equation~\eqref{eq:road-filter-predict} then propagates this posterior to $b_{t+1}^\Theta$, the distribution of $\Theta_{t+1}$ that the planner uses at decision epoch $t+1$. Online Supplement Section~\ref{app:road-filter} derives the update--predict recursion. Because the road state takes finitely many values, the planner can track road conditions exactly. The computational difficulty of the problem lies elsewhere, namely, in the continuous compliance belief and the variable-dimension traffic state.

\subsection{Traffic Dynamics, Performance Feedback, and Cost}
\label{subsec:dp-traffic-cost}

We now describe how the realized route choices move through the network. Newly departing travelers enter the first arc of their selected route during period $t$. Using the arc identifier $e_{p,k}$ from Section~\ref{subsec:dp-state}, we combine new departures with inherited entries to obtain
\begin{equation}\label{eq:dp-arc-flow}
   z_{ij,t}
   = h_{ij,t}(\mathbf X_t)
   +
    \sum_{o\in\mathcal O_t^+(S_t)}
    \sum_{p\in\mathcal P_{o,t}(S_t)}
    \left[
        Y_{op,t}\mathbf 1\{e_{p,1}=(i,j)\}
        +
        (N_{op,t}-Y_{op,t})\mathbf 1\{e_{p^0_{op,t},1}=(i,j)\}
    \right],
\end{equation}
for all $(i,j)\in\mathcal E$. In road state $\vartheta$, the arc travel-time function is
\begin{equation}
\label{eq:dp-travel-time}
    \operatorname{tt}^{\vartheta}_{ij,t}(z)
    =
    t^0_{ij,t}(\vartheta)
    \left[
        1+\alpha_{\mathrm{BPR}}
        \left(\frac{z}{c_{ij,t}(\vartheta)}\right)^{
            \beta_{\mathrm{BPR}}}
    \right].
\end{equation}
Here, $t^0_{ij,t}(\vartheta)>0$ is the free-flow travel time, $c_{ij,t}(\vartheta)>0$ is capacity expressed as a traveler count per decision period, and $\alpha_{\mathrm{BPR}}\geq0$ and $\beta_{\mathrm{BPR}}>0$ are the Bureau of Public Roads (BPR) parameters~\citep{manual1964bureau}. For period length $\Delta t>0$ in the same time units, an entering packet receives the fixed traversal duration
\begin{equation}
\label{eq:packet-duration}
    \ell^{\mathrm{trav}}_{ij,t}
    =
    \left\lceil
        \frac{\operatorname{tt}^{\Theta_t}_{ij,t}(z_{ij,t})}
             {\Delta t}
    \right\rceil.
\end{equation}
For each recommendation $(o,p)$, we create a packet of $Y_{op,t}$ travelers on route $p$ and a packet of $N_{op,t}-Y_{op,t}$ travelers on route $p^0_{op,t}$, omitting packets of size zero. Each new packet enters its first arc. An inherited packet with zero residual duration enters its next arc. We compute the combined arc-entry counts using~\eqref{eq:dp-arc-flow} and assign each entering packet the duration in~\eqref{eq:packet-duration}. Its residual duration at the next decision epoch is $\ell^{\mathrm{trav}}_{ij,t}-1$. An inherited packet with a positive residual duration stays on its current arc and decreases that residual by one. We remove packets that complete their final arc and form $\mathbf X_{t+1}$ from the packet records that result from these operations.

We collect the physical and behavioral outcomes of period $t$, including the current hidden road state, in
\begin{equation}
\label{eq:random-information}
    \Xi_{t+1}
    =
    \left(
        \mathbf D_t^{\mathrm{act}},
        \mathbf D_{t+1}^{\mathrm{new}},
        \mathbf N_t,
        \mathbf Y_t,
        \Theta_t
    \right).
\end{equation}
The stage cost charges unrounded travel time at each arc entry:
\begin{equation}
\label{eq:dp-stage-cost}
    g_t(S_t,\mathbf r_t,\Xi_{t+1})
    =
    \sum_{(i,j)\in\mathcal E}
    z_{ij,t}\operatorname{tt}^{\Theta_t}_{ij,t}(z_{ij,t}), 
\end{equation}
which serves as a travel-time proxy based on the unrounded entry-time function, whereas physical propagation uses the rounded durations in~\eqref{eq:packet-duration}. We charge an arc entered before $T$ in full even if its traversal ends after $T$. The baseline model includes no additional period cost for origin waiting or incident exposure.

The performance state closes the feedback loop between the quality of guidance and future compliance. To evaluate guidance, let $\widehat{\operatorname{tt}}^{\mathrm{snap}}_{ij,t}$ be the observed or reconstructed period-$t$ arc-time snapshot in $O^\Theta_{t+1}$, and define $\widehat{\operatorname{TT}}^{\mathrm{snap}}_{p,t} =\sum_{(i,j)\in p}\widehat{\operatorname{tt}}^{\mathrm{snap}}_{ij,t}$. We evaluate the recommended and default routes on this same snapshot, which gives a contemporaneous comparison for trips that span multiple periods. 
Let $N^{\mathrm{inf}}_{o,t} =\sum_{p\in\mathcal P^{\mathrm{inf}}_{o,t}(S_t)}N_{op,t}$ for active groups and $N^{\mathrm{inf}}_{o,t}=0$ otherwise. For relative tolerance $\varepsilon_A\geq0$ and $N^{\mathrm{inf}}_{o,t}>0$, we define
\begin{equation}
\label{eq:performance-signal}
    a^{\mathrm{perf}}_{o,t}
    =
    \frac{1}{N^{\mathrm{inf}}_{o,t}}
    \sum_{p\in\mathcal P^{\mathrm{inf}}_{o,t}(S_t)}
    N_{op,t}
    \mathbf 1\!\left\{
        \widehat{\operatorname{TT}}^{\mathrm{snap}}_{p,t}
        \leq
        (1+\varepsilon_A)
        \widehat{\operatorname{TT}}^{\mathrm{snap}}_{p^0_{op,t},t}
    \right\}.
\end{equation}
We exclude recommendations whose route coincides with the default, because they provide no assessment of nontrivial guidance. The signal $a^{\mathrm{perf}}_{o,t}$ is therefore the fraction of group $o$'s informative recommendations that performed within the tolerance of their default routes. With smoothing weight $\alpha_\tau\in[0,1]$, the recommendation performance state evolves following
\begin{equation}
\label{eq:performance-ema}
    \tau_{o,t+1}
    =
    \begin{cases}
        (1-\alpha_\tau)\tau_{o,t}
        +\alpha_\tau a^{\mathrm{perf}}_{o,t},
        &N^{\mathrm{inf}}_{o,t}>0,\\
        \tau_{o,t},
        &N^{\mathrm{inf}}_{o,t}=0.
    \end{cases}
\end{equation}
The update preserves $\tau_{o,t}\in[0,1]$ and leaves it unchanged when the planner issues no informative recommendations to group $o$, which includes the case of an inactive group. Each group shares this performance score across its acceptable destinations. Because $\tau_{o,t}$ enters the compliance features in~\eqref{eq:compliance-feature-vector}, a run of recommendations that underperform their defaults lowers future compliance whenever travelers reward reliable guidance, that is, whenever the performance coefficient is positive. A run of reliable recommendations raises it.

Finally, for penalty weights $\eta_q,\eta_X\geq0$, the terminal penalty is
\begin{equation}
\label{eq:dp-terminal-cost}
    G_T(S_T)
    =
    \eta_q\sum_{o\in\mathcal O}q_{o,T}
    +\eta_X U(\mathbf X_T),
    \quad
    U(\mathbf X_T)
    =
    \sum_{c=1}^{C_T}n_{c,T}.
\end{equation}
The weights convert waiting and unfinished traveler counts into stage-cost units. Under a common exogenous demand realization, the queue penalty is identical across recommendation policies, whereas the unfinished-traffic penalty depends on route choices. The latter penalty prevents the planner from shifting travel cost beyond the planning horizon.

\subsection{Controlled Transition and Bellman Recursion}
\label{subsec:dp-transition}

We now assemble the pieces. Let $\mathbf O_{t+1}  = (O^C_{t+1},O^\Theta_{t+1})$ collect the post-period observations used for compliance learning and road state filtering. The origin queues satisfy
\begin{equation}\label{eq:q-transition}
    q_{o,t+1}=q_{o,t}-D^{\mathrm{act}}_{o,t}+D^{\mathrm{new}}_{o,t+1},
    \quad o\in\mathcal O.
\end{equation}
The packet propagation in Section~\ref{subsec:dp-traffic-cost} produces $\mathbf X_{t+1}$ using~\eqref{eq:dp-arc-flow}--\eqref{eq:packet-duration}, and Equations~\eqref{eq:bayes-compliance-update}, \eqref{eq:road-filter-update}--\eqref{eq:road-filter-predict}, and~\eqref{eq:performance-ema} update $b_{t+1}^C$, $b_{t+1}^\Theta$, and $\boldsymbol\tau_{t+1}$, respectively. Time then advances to $t+1$. A deterministic function $f_t$ assembles these components from the current state, action, outcomes, and observations, and we write $S_{t+1}=f_t(S_t,\mathbf r_t,\Xi_{t+1},\mathbf O_{t+1})$.

A planner who must track a joint belief over the compliance coefficients and the road state faces a much larger information state than one who can track the two separately. Proposition~\ref{prop:belief-factorization} shows that, under our assumptions, separate beliefs suffice.

\begin{proposition}[Conditional Belief Factorization]
\label{prop:belief-factorization}
Under the initial independence assumption~\eqref{eq:initial-belief-factorization} and the transition and observation assumptions in Section~\ref{sec:dynamic-route-model}, the compliance coefficients $\boldsymbol\beta^C$ and current road state $\Theta_t$ remain conditionally independent given the observable history $\mathscr F_t$. At every decision epoch $t$, for every measurable $B\subseteq\mathbb R^5$ and $\vartheta\in\mathcal S^\Theta$,
\begin{equation}
\label{eq:belief-factorization}
    \Pr\!\left(
        \boldsymbol\beta^C\in B,\Theta_t=\vartheta
        \mid\mathscr F_t
    \right)
    =
    b_t^C(B)b_t^\Theta(\vartheta).
\end{equation}
\end{proposition}

Online Supplement Section~\ref{app:belief-factorization-proof} proves Proposition~\ref{prop:belief-factorization}. The factorization lets the planner learn about travelers and about roads through two separate and simpler updates.

Let $P_t$ denote the joint conditional distribution of $(\Xi_{t+1},\mathbf O_{t+1})$. Integrating over the compliance and road state beliefs gives the predictive distribution
\begin{equation}
\label{eq:predictive-distribution}
    P_t(d\zeta,d\boldsymbol\xi\mid S_t,\mathbf r_t)
    =
    \int_{\mathbb R^5}
    \sum_{\vartheta\in\mathcal S^\Theta}
    P_t(
        d\zeta,d\boldsymbol\xi
        \mid
        S_t,\mathbf r_t,
        \boldsymbol\beta^C,
        \Theta_t=\vartheta
    )
    b_t^C(d\boldsymbol\beta^C)
    b_t^\Theta(\vartheta).
\end{equation}
Here $\zeta$ and $\boldsymbol\xi$ denote realizations of $\Xi_{t+1}$ and $\mathbf O_{t+1}$, respectively. The conditional distribution on the right-hand side combines the demand distribution in Section~\ref{subsec:dp-demand-counts}, the assignment and agreement distributions in~\eqref{eq:recommendation-multinomial} and~\eqref{eq:compliance-binomial}--\eqref{eq:coincident-route-agreement}, and the road observation distribution in~\eqref{eq:road-observation-kernel}. Equation~\eqref{eq:predictive-distribution} averages the current road state $\Theta_t$ over its belief, and the transition kernel~\eqref{eq:road-transition-kernel} enters the next-state map through the belief prediction~\eqref{eq:road-filter-predict}. Online Supplement Section~\ref{app:road-filter} discusses a case in which compliance coefficient and road state are conditionally dependent.

We assume that the route sets are finite, the action correspondence is measurable with nonempty compact values, the transition and observation kernels are measurable, expected costs are finite, and the standard lower semicontinuity and attainment conditions for finite-horizon controlled Markov models hold.

\begin{proposition}[Belief-State Sufficiency]
\label{prop:belief-state-sufficiency}
Under the transition and observation models and the conditional independence assumptions in Section~\ref{sec:dynamic-route-model}, the state $S_t$ in~\eqref{eq:dp-state} is a controlled Markov information state. In particular, for every admissible period-$t$ action $\mathbf r_t\in\mathcal A(S_t)$ and all measurable sets $\mathcal B$ and $\mathcal C$,
\begin{equation}\label{eq:controlled-markov-property}
    \Pr(S_{t+1}\in\mathcal B,\ g_t\in\mathcal C\mid \mathscr F_t,\mathbf r_t)
    =
    \Pr(S_{t+1}\in\mathcal B,\ g_t\in\mathcal C\mid S_t,\mathbf r_t)
    \quad\text{almost surely}.
\end{equation}
Consequently, restricting attention to Markov policies entails no loss of optimality. Under the stated regularity conditions, an optimal measurable Markov selector exists for the exact information state $S_t$.
\end{proposition}
Online Supplement Section~\ref{app:belief-state-proof} proves Proposition~\ref{prop:belief-state-sufficiency}. The proposition has a practical reading. The planner need not remember the full history of recommendations, departures, and route choices. Waiting demand, the traffic already in the network, two beliefs, and a performance score carry everything that matters for future decisions.

Let $V_t^*(S_t)$ denote the minimum expected cost from period $t$ onward, with terminal value $V_T^*(S_T)=G_T(S_T)$. For $t=T-1,\ldots,0$, we define the optimal action-value function
\begin{equation}
\label{eq:optimal-q-function}
    Q_t^*(S_t,\mathbf r_t)
    =
    \mathbb E\!\left[
        g_t(S_t,\mathbf r_t,\Xi_{t+1})
        +
        V_{t+1}^*(S_{t+1})
        \mid S_t,\mathbf r_t
    \right],
\end{equation}
where the expectation uses the predictive distribution in~\eqref{eq:predictive-distribution}. The Bellman recursion and an optimal policy are
\begin{align}
\label{eq:bellman-main}
    V_t^*(S_t)
    & =
    \min_{\mathbf r_t\in\mathcal A(S_t)}
    Q_t^*(S_t,\mathbf r_t),\\
\label{eq:optimal-policy}
    \mu_t^*(S_t)
    & \in
    \arg\min_{\mathbf r_t\in\mathcal A(S_t)}
    Q_t^*(S_t,\mathbf r_t).
\end{align}
The recursion captures three ways in which a recommendation matters: its current arc-entry cost, its effect on future traffic propagation, and its effect on what the planner learns and on how travelers respond to later guidance. A myopic rule captures only the first. Solving this recursion exactly, however, is computationally expensive because the packet state has variable dimension, the compliance belief is distribution-valued, and the route-share action lies in a continuous simplex. In Section~\ref{sec:ppo-approximation}, we therefore develop a particle-based PPO approximation. 

\section{Approximate Solution Method}
\label{sec:ppo-approximation}

\begin{figure}[htbp!]
  \centering
  \includegraphics[width=0.8\linewidth]{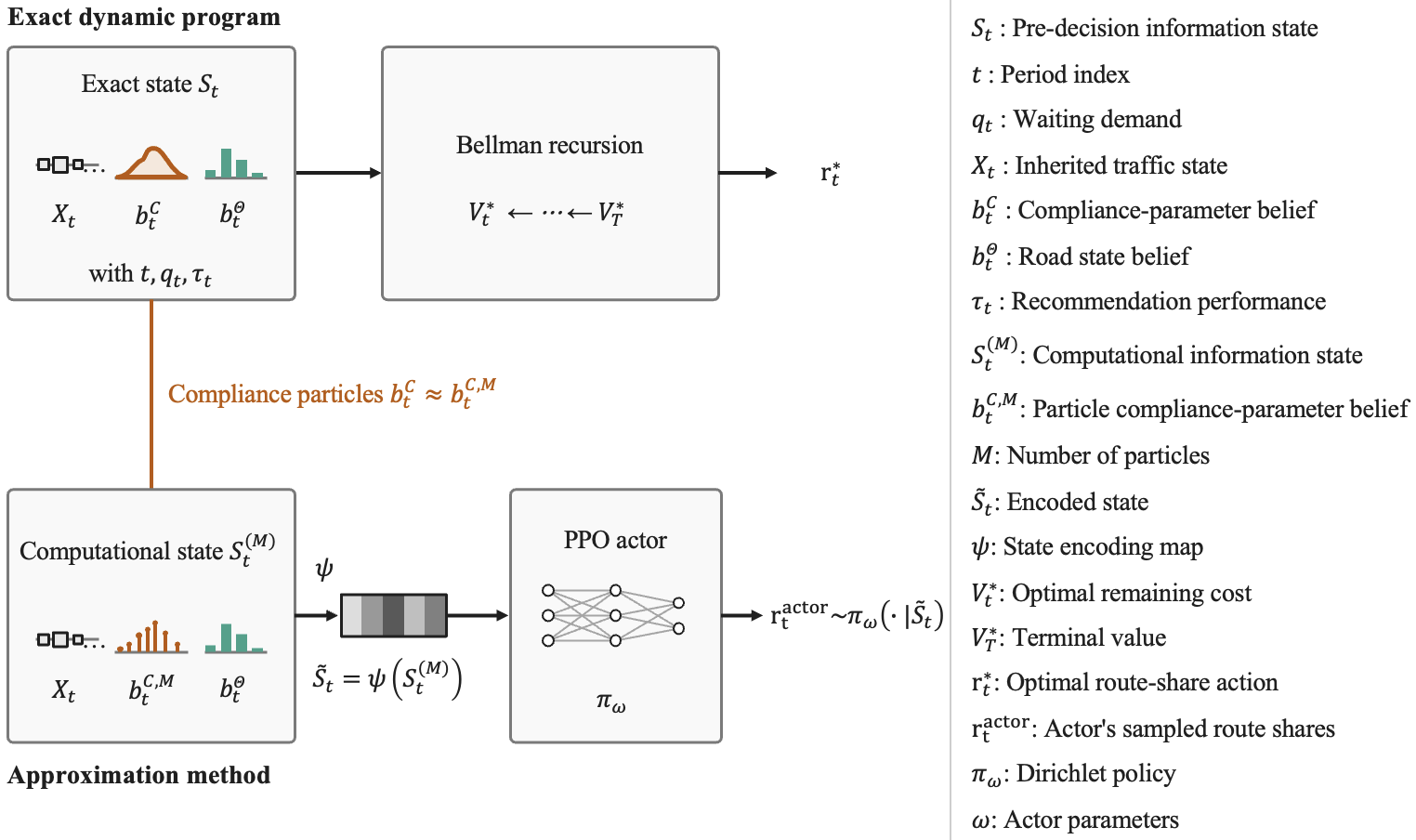}
  \caption{Exact Bellman policy and computational approximation. We replace the compliance posterior with weighted particles and retain the remaining state components. The Dirichlet actor uses a fixed-dimensional state encoding to sample route shares $\mathbf r_t^{\mathrm{actor}}$ during training. The planner forms the feasible action $\mathbf r_t$ from these shares as described in Section~\ref{subsec:ppo-actor}.}
  \label{fig:route-guidance-approximation}
\end{figure}

The Bellman recursion specifies what the planner should optimize, but its three sources of complexity prevent exact computation on networks of practical size. We keep the model of Section~\ref{sec:dynamic-route-model} intact and approximate each source of complexity in turn. Specifically, we approximate the Bellman policy in~\eqref{eq:optimal-policy} with weighted particles for the posterior of the continuous compliance coefficients (Section~\ref{subsec:particle-approximation}), a fixed-dimensional state encoding (Section~\ref{subsec:state-encoding}), and a Dirichlet actor and value-function critic trained on simulated trajectories with proximal policy optimization (PPO; \citealp{schulman2017proximal}; Sections~\ref{subsec:ppo-actor}--\ref{subsec:ppo-clipped-objective}). These approximations represent beliefs, policies, and values, and they compute actions for the model in Section~\ref{sec:dynamic-route-model}. Figure~\ref{fig:route-guidance-approximation} contrasts the Bellman recursion with these approximations.

\subsection{Particle Approximation of the Compliance Belief}
\label{subsec:particle-approximation}

For an integer $M\geq 1$, we approximate the posterior in~\eqref{eq:bayes-compliance-update} with the empirical measure
\begin{equation}
\label{eq:particle-compliance-belief}
    b_t^{C,M}(d\boldsymbol\beta^C)
    =
    \sum_{m=1}^{M}w_t^{(m)}
    \delta_{\boldsymbol\beta_t^{C,(m)}}(d\boldsymbol\beta^C),
    \quad
    w_t^{(m)}\geq0,
    \quad
    \sum_{m=1}^{M}w_t^{(m)}=1,
\end{equation}
where $\delta_x$ denotes a point mass at $x$ and $\boldsymbol\beta_t^{C,(m)}\in\mathbb R^5$ denotes particle $m$ at period $t$. We draw the initial particles independently from $b_0^C$ and assign them equal weights $1/M$. Each particle represents one hypothesis about how travelers respond to guidance, and its weight records how well that hypothesis has explained the route choices observed so far. Replacing $b_t^C$ with $b_t^{C,M}$ in~\eqref{eq:dp-state} gives the computational state
$S_t^{(M)}=(t,\mathbf q_t,\mathbf X_t,b_t^{C,M},b_t^\Theta,\boldsymbol\tau_t)$.

For the particle update, the likelihood~\eqref{eq:compliance-likelihood} uses pre-decision features from $S_t^{(M)}$ instead of $S_t$, and the resulting likelihood $L_t^{C,M}$ holds these features fixed across candidate compliance coefficients. After observing $O^C_{t+1}$, we evaluate $L_t^{C,M}$ at each particle using the integer counts $(\mathbf N_t,\mathbf Y_t)$ and normalize the updated weights:
\begin{equation}
\label{eq:particle-weight-update}
    \bar w_{t+1}^{(m)}
    =
    \frac{
        w_t^{(m)}
        L_t^{C,M}(\boldsymbol\beta_t^{C,(m)})
    }{
        \displaystyle\sum_{j=1}^{M}
        w_t^{(j)}
        L_t^{C,M}(\boldsymbol\beta_t^{C,(j)})
    },
    \qquad m=1,\ldots,M.
\end{equation}

The effective sample size $\operatorname{ESS}_{t+1} = 1/ \sum_{m=1}^{M}\left(\bar w_{t+1}^{(m)}\right)^2 $ measures how many particles still carry meaningful weight. If $\operatorname{ESS}_{t+1}<M_{\mathrm{ESS}}$, where $M_{\mathrm{ESS}}\in[0,M]$ is the resampling threshold, we resample $M$ particles from the updated empirical measure and reset their weights to $1/M$. Otherwise, we retain the particles and set $w_{t+1}^{(m)}=\bar w_{t+1}^{(m)}$. 
We store $b_t^\Theta$ as a finite probability vector and update it through~\eqref{eq:road-filter-update}--\eqref{eq:road-filter-predict} using the realized observation context, without a particle approximation.

    \subsection{Fixed-Dimensional State Encoding}
\label{subsec:state-encoding}

The packet state and the particle cloud vary in size, whereas a neural network requires an input of fixed length. In PPO, the actor and critic therefore receive the fixed-dimensional vector $\widetilde S_t = \psi(S_t^{(M)}) \in\mathbb R^{d_\psi}$, where a specified deterministic function $\psi$ extracts computational-state features and $d_\psi$ denotes its output dimension. For instance, $\psi$ may combine observed route choices, traffic and belief summaries, route predictions, and route-availability indicators. Because $\psi$ need not preserve all information in $S_t^{(M)}$, three approximation errors may arise relative to the Bellman policy: (i) finite-particle error in $b_t^{C}$, (ii) information loss through $\psi$, and (iii) policy and value-function approximation and optimization error. Experiment E1 in Section~\ref{subsec:exp-e1} assesses the combined effect of approximation and optimization error against an exact benchmark.

\subsection{Dirichlet Policy Parameterization}
\label{subsec:ppo-actor}

A route-share action is a probability vector for each group, and thus the actor needs a distribution over the simplex. The Dirichlet distribution fits this need because every sample is a feasible share vector by construction. The actor samples route shares for every group with at least two available routes, which form the set
$\mathcal O_t^{\geq2}(S_t^{(M)})=\{o\in\mathcal O:|\mathcal P_{o,t}(S_t^{(M)})|\geq2\}$.
For $o\in\mathcal O_t^{\geq2}(S_t^{(M)})$ and $p\in\mathcal P_{o,t}(S_t^{(M)})$, the actor with parameters $\boldsymbol\omega$ maps its unconstrained output $y_{op,t}(\widetilde S_t;\boldsymbol\omega)\in\mathbb R$ into a positive Dirichlet concentration parameter \citep{ng2011dirichlet}:
\begin{equation}
\label{eq:dirichlet-concentration}
    \kappa_{op,t}(\widetilde S_t;\boldsymbol\omega)
    =
    \log\!\left(
        1+e^{y_{op,t}(\widetilde S_t;\boldsymbol\omega)}
    \right)
    +\kappa_{\min},
    \quad \kappa_{\min}>0.
\end{equation}
Let $\boldsymbol\kappa_{o,t}$ collect these parameters. Their normalized values determine the mean shares, and for fixed mean shares, a larger sum of concentrations reduces variability \citep{minka2000estimating}. During training, the actor therefore learns both where to send travelers and how much to explore around that allocation. Conditional on the feature vector and the available route sets from $S_t^{(M)}$, the actor samples shares independently across groups:
\begin{equation}
\label{eq:dirichlet-policy-origin}
    \mathbf r_{o,t}^{\mathrm{actor}}
    \sim
    \operatorname{Dirichlet}\!\left(
        \boldsymbol\kappa_{o,t}(\widetilde S_t;\boldsymbol\omega)
    \right),
    \quad o\in\mathcal O_t^{\geq2}(S_t^{(M)}).
\end{equation}
Let $\mathbf r_t^{\mathrm{actor}}$ collect the sampled group vectors, with $\mathbf r_{o,t}^{\mathrm{actor}}=(r_{op,t}^{\mathrm{actor}})_{p\in\mathcal P_{o,t}(S_t^{(M)})}$. Suppressing the dependence on the available route sets, we write the joint sampling density as
\begin{equation}
\label{eq:joint-dirichlet-policy}
    \pi_{\boldsymbol\omega}
    (\mathbf r_t^{\mathrm{actor}}\mid\widetilde S_t)
    =
    \prod_{o\in\mathcal O_t^{\geq2}(S_t^{(M)})}
    \operatorname{Dirichlet}\!\left(
        \mathbf r_{o,t}^{\mathrm{actor}};
        \boldsymbol\kappa_{o,t}(\widetilde S_t;\boldsymbol\omega)
    \right).
\end{equation}
If the index set is empty, the actor draws no shares, and the empty product equals one.

The planner forms the feasible action $\mathbf r_t\in\mathcal A(S_t^{(M)})$ from the actor's shares. For each active group $o\in\mathcal O_t^+(S_t^{(M)})$, it sets $\mathbf r_{o,t} = \mathbf r_{o,t}^{\mathrm{actor}}$ if the group has at least two available routes and $\mathbf r_{o,t} = (1)$ if the group has one available route. Groups with empty queues have no entry in $\mathbf r_t$. The simulator uses $\mathbf r_t$ to assign departing travelers through~\eqref{eq:recommendation-multinomial}. PPO uses the full sampled vector $\mathbf r_t^{\mathrm{actor}}$ in the likelihood ratio and includes every Dirichlet factor in the entropy term (Section~\ref{subsec:ppo-clipped-objective}).

For deterministic policy evaluation, we replace the sampled shares with their Dirichlet means:
\begin{equation}
\label{eq:dirichlet-mean-action}
    r_{op,t}^{\mathrm{actor,eval}}
    =
    \frac{\kappa_{op,t}}
    {\sum_{p'\in\mathcal P_{o,t}(S_t^{(M)})}
        \kappa_{op',t}},
    \quad
    o\in\mathcal O_t^{\geq2}(S_t^{(M)}),
    \quad
    p\in\mathcal P_{o,t}(S_t^{(M)}).
\end{equation}
We form the action from these mean shares using the same rule and evaluate it under the stochastic traffic model. For each group with a Dirichlet factor, sampled shares are almost surely positive on every available route, and mean shares are strictly positive, which excludes exact boundary allocations. Because every nondefault route receives a positive share, the policy issues informative recommendations with positive probability whenever travelers depart, and thus it keeps learning about compliance throughout the horizon.

\subsection{Reward, Training, and Policy Execution}
\label{subsec:ppo-reward-reformulation}
\label{subsec:ppo-clipped-objective}

The critic $V_{\boldsymbol\nu}(\widetilde S_t)$, with parameters $\boldsymbol\nu$, estimates the current policy's expected remaining return. In our implementation, the actor and critic share hidden layers and have separate output heads. Thus, $\boldsymbol\omega$ and $\boldsymbol\nu$ both include the shared hidden-layer parameters, and each includes its own output-head parameters. For $\gamma=1$ and untransformed rewards, the true value of the current policy at $t<T$ equals the negative expected remaining cost, including the terminal penalty.
(Online Supplement Section~\ref{app:ema-implementation} specifies the reward transformations used in the experiments.) 

We define the reward and terminal critic value by
\begin{equation}
\label{eq:ppo-reward}
    R_{t+1}
    =
    -g_t(S_t^{(M)},\mathbf r_t,\Xi_{t+1})
    -
    \mathbf 1\{t=T-1\}G_T(S_T^{(M)}),
    \quad
    V_{\boldsymbol\nu}(\widetilde S_T)=0.
\end{equation}
We seek actor parameters that maximize the expected return
\begin{equation}
\label{eq:rl-objective}
    \max_{\boldsymbol\omega}\;
     \mathbb E^{\pi_{\boldsymbol\omega}}
    \left[
        \sum_{t=0}^{T-1}\gamma^tR_{t+1}
    \right],
\end{equation}
where $\gamma$ is the discount factor. We set $\gamma=1$ to match the undiscounted finite-horizon objective in~\eqref{eq:problem-objective}. Summing~\eqref{eq:ppo-reward} over the horizon gives
\begin{equation}
\label{eq:reward-cost-identity}
    \sum_{t=0}^{T-1}R_{t+1}
    =
    -\left[
        \sum_{t=0}^{T-1}g_t(S_t^{(M)},\mathbf r_t,\Xi_{t+1})
        +G_T(S_T^{(M)})
    \right].
\end{equation}
Within a fixed policy class, maximizing the unregularized expected return is therefore equivalent to minimizing the finite-horizon expected cost. In particular, the learned policy faces the same trade-off between current travel cost and unfinished traffic at time $T$ as the exact model.

We train the actor and critic with the clipped PPO objective of \citet{schulman2017proximal} and generalized advantage estimates. Within each rollout batch, a critic held at its batch-start parameters computes the advantages and critic targets. The likelihood ratio uses the full sampled share vector $\mathbf r_t^{\mathrm{actor}}$, and the entropy bonus sums the differential entropies of the Dirichlet factors in~\eqref{eq:joint-dirichlet-policy}. Clipping keeps each update close to the policy that generated the batch, which stabilizes training when simulated costs are noisy. Online Supplement Section~\ref{app:ppo-training} states the advantage estimator, the clipped objective, and the batch protocol.

We train the actor and critic offline on simulated trajectory batches and fix their parameters during validation and execution. Within each trajectory, we update particle weights through~\eqref{eq:particle-weight-update}, road beliefs through~\eqref{eq:road-filter-update}--\eqref{eq:road-filter-predict}, packet traffic through Section~\ref{subsec:dp-traffic-cost}, origin queues through~\eqref{eq:q-transition}, and recommendation performance through~\eqref{eq:performance-ema}. The performance state then enters the next period's compliance feature vector in~\eqref{eq:compliance-feature-vector}. This division of labor places the computational burden before an emergency occurs. During an operation, the planner only updates its beliefs and evaluates the trained network at each decision epoch.

\section{Numerical Experiments}
\label{sec:numerical}

We design the experiments around three questions that a planner would ask before adopting our approach. First, does PPO find near-optimal recommendations when we can compute the optimum exactly? Second, on a realistic network, does planning over the full horizon outperform simpler rules, including rules that know the true compliance behavior? Third, how much of the benefit comes from updating beliefs and recommendation performance online, and does the answer depend on the quality of the initial prior? Experiments E1, E2, and E3 address these questions, respectively. Section~\ref{subsec:exp-setup} introduces the Eastern Massachusetts (EMA) instance and metrics for evaluation. In E1, we compare a variant of PPO that uses a categorical distribution over nine discrete route-share allocations with an exact dynamic program over the same action set (Section~\ref{subsec:exp-e1}). In E2, we compare the mean and upper-tail costs of PPO with those of benchmark policies (Section~\ref{subsec:exp-e2}). In E3, we use the same Dirichlet policy parameterization to compare three information settings under three initial compliance priors (Section~\ref{subsec:exp-e3}). Online Supplement Section~\ref{app:terminal-sensitivity} assesses sensitivity to the terminal penalty.

\subsection{EMA Instance and Statistical Protocol}
\label{subsec:exp-setup}

We use the EMA highway data in the Transportation Networks repository \citep{codetransportation} for E2 and E3. The source network contains $74$ nodes and $258$ directed links. We select four origin--destination (OD) groups with fixed destinations 
and three candidate routes per group. These twelve routes use $13$ distinct links and define the bottleneck instance for our evaluation. Each group initially has $1,053$ waiting travelers, the network starts empty, and the horizon comprises $T=12$ periods of $15$ minutes each. (Online Supplement Section~\ref{app:ema-implementation} details the period-specific binomial departures and Poisson arrivals to origin queues.) In E2 and E3, all three routes per group remain available throughout the horizon, and therefore the joint density~\eqref{eq:joint-dirichlet-policy} includes all four groups at every decision epoch.

The instance represents a pattern that makes evacuation guidance valuable: travelers' habitual routes converge on one critical link. All four default routes share one arc. 
An incident multiplies that arc's free-flow time by $3.5$ and its capacity by $0.2$. We set the initial incident probability to $0.30$ and the normal-to-incident and incident-to-incident transition probabilities to $0.18$ and $0.75$, respectively. Incidents change travel times and capacities without closing routes. Each group's first bypass has a free-flow travel time $0.43$--$0.52\%$ greater than its default route. The second alternative has a longer free-flow travel time than the default and still uses the incident-sensitive arc for two groups. This setting provides a bypass with little additional free-flow travel time and tests how the planner allocates recommendations across routes over time under congestion, uncertain road conditions, and partial compliance. All links use the BPR parameters $(\alpha_{\mathrm{BPR}},\beta_{\mathrm{BPR}})=(0.15,4)$.

In E2 and E3, we fix the true compliance vector at $\boldsymbol\beta^C=(-0.20,-2.00,1.10,1.00,-0.70)^\top$ throughout training and evaluation, and thus the initial priors are alternative beliefs about identical coefficients. Under these coefficients, travelers are less likely to comply when a recommendation is predicted to be slower or more congested than the alternative, and more likely to comply when they know the recommended route or when recent guidance has performed well. Updating policies represent compliance beliefs with $3000$ particles and resample systematically when the effective sample size falls below $1500$. We set the performance parameters to $\tau_{o,0}=0.60$, $\alpha_\tau=0.20$, and $\varepsilon_A=0.05$, and the feature parameters to $s_{\mathrm{TT}}=0.25$ hours and $C_{\mathrm{Cong}}=5$. Although the response coefficients are fixed, compliance probabilities vary with traffic and recommendation performance.
%
%
After each period, updating policies receive an incident alarm with sensitivity $0.85$ and false-positive probability $0.10$, together with an error-free link-time snapshot that identifies the completed period's road state given realized flows. They update road beliefs through~\eqref{eq:road-filter-update} and predict the road state at the next decision epoch through~\eqref{eq:road-filter-predict}. A current-road oracle observes the state before acting.

In E2 and E3, our primary performance measure is total cost in hours, which is the sum of arc-entry travel costs~\eqref{eq:dp-stage-cost} and the terminal penalty $G_T(S_T)=20\sum_{o\in\mathcal O}q_{o,T}+5U(\mathbf X_T)$. Each arc entry contributes its entering traveler count multiplied by the arc's unrounded travel time. At the end of the three-hour horizon, the terminal penalty adds $20$ hours per traveler still waiting at an origin and $5$ hours per departed traveler still en route. We report this cost before reward normalization or clipping. Because the cost includes these penalties, it does not measure the elapsed time until the last traveler reaches a destination. We report terminal traveler counts separately when we explain cost differences. 

We compare PPO with nine benchmark policies in E2. These benchmarks represent rules that a planner could deploy without a learned policy. Five use route-aware myopic optimization: (i) static myopic, (ii) frozen-information myopic, (iii) updating myopic, (iv) oracle-compliance myopic, and (v) oracle-compliance-and-road myopic. These policies use coordinate descent over route shares in increments of $0.2$ to minimize a full-route cost proxy for current departures. The proxy uses inherited arc entries and expected departures after compliance to price every arc on each route. The static policy (i) repeats its initial allocation. The frozen-information policy (ii) recomputes allocations using the initial beliefs and recommendation performance. The updating policy (iii) applies the compliance update~\eqref{eq:bayes-compliance-update} through its particle implementation~\eqref{eq:particle-weight-update}, the road updates~\eqref{eq:road-filter-update}--\eqref{eq:road-filter-predict}, and the performance update~\eqref{eq:performance-ema} after each period. The two oracle variants (iv) and (v) additionally receive either the true compliance vector or both that vector and the current road state. Four further benchmarks recommend only default routes, allocate shares uniformly, sample random Dirichlet shares, or allocate shares by a softmax over negative predicted route travel time with a temperature of $0.40$ hours. The default-only policy reproduces the outcome without guidance, because every traveler then takes the default route. Online Supplement Section~\ref{app:ema-evaluation} defines the myopic proxy and the softmax rule. Queues and traffic evolve under every policy. Our comparisons concern these specific implementations, which differ in planning horizon, action discretization, and optimization procedure.

For PPO in E2 and E3, we validate and evaluate the Dirichlet mean policy in~\eqref{eq:dirichlet-mean-action}. We use out-of-sample simulation scenarios that differ from those used for training and checkpoint selection, compute each actor's performance metrics separately, and then average the actor-specific estimates. The metrics include mean total cost, the empirical cost quantiles $\widehat q_{0.90}$ and $\widehat q_{0.99}$ (the $90$th and $99$th percentiles), and the average total cost among the highest-cost $5\%$ of scenarios, denoted by $\widehat{\CVaR}_{0.95}$. The tail metrics matter to an emergency manager because a policy that performs well on average but poorly in the worst scenarios fails when the stakes are highest. 
Let $\Delta_a$ denote the scenario-average paired cost difference for training replicate $a=1,2,3$. We specify the subtraction order for each comparison. We compute the $95\%$ Student-$t$ interval across training replicates as $\bar\Delta \pm t_{0.975,2}(s_\Delta/\sqrt{3})$, where $\bar\Delta$ and $s_\Delta$ are the sample mean and standard deviation of the three differences, and $t_{0.975,2}$ is the $97.5$th percentile of the Student-$t$ distribution with two degrees of freedom. In E2, these intervals summarize variation across training replicates conditional on a common evaluation set. In E3, they also reflect variation across the distinct evaluation sets for each replicate. 

\subsection{E1: Exact Computation and Validation}
\label{subsec:exp-e1}

E1 asks whether PPO finds the best recommendations in a setting where we can verify the answer. Specifically, it tests whether categorical PPO can recover a minimum-cost schedule from a finite set of route-share choices in a deterministic instance. The instance has two origins and two routes per origin, with departures $(8,2,5)$ from $O_1$ and $(2,8,5)$ from $O_2$ over three decision periods. Because the two origins' departures peak in different periods, a good schedule may need to change over time. We convert shares to integer recommendations by largest-remainder apportionment, require all travelers to follow their assigned routes, and hold the compliance belief, road belief, and recommendation performance state fixed. The road stays in its normal state. We replace the baseline terminal penalty with the arc-entry cost of a $15$-period clearance continuation with no new departures. 

Each origin's default-route share lies in $\{0,0.5,1\}$, which yields nine joint route-share choices per period. The Bellman recursion gives the same minimum total travel cost as an exhaustive evaluation of all $9^3=729$ three-period schedules. We compare all nine repeated-action policies, a one-period policy, and five independently trained categorical PPO actors. PPO samples these nine actions during training and selects the most probable action during evaluation. Online Supplement Section~\ref{app:e1-details} details calibration, training, checkpoint selection, and actor-level results.

Table~\ref{tab:exp-e1-gap} compares the exact dynamic program (DP) with the PPO-based method, the best static route allocation, and the one-period policy.
The PPO row averages the selected checkpoint costs of the five independently trained actors.

\begin{table}[htbp!]
  \centering
   \caption{E1 validation over nine route-share choices per period, reporting ``Cost'' (in hours), the cost difference from the minimum given by exact DP (``Cost minus minimum''), and the number of PPO actors that recover the exact minimum cost (``Exact recovery''). 
   }
   \label{tab:exp-e1-gap}
  \small
  \begin{tabular}{lrrc}
    \toprule
    Policy & Cost & Cost minus minimum & Exact recovery \\
    \midrule
    \emph{Exact discrete-action DP} & $121.60$ & --- & --- \\
    \emph{Categorical PPO}       & $121.75$ & $0.15$ & $4/5$ \\
    \emph{Best static allocation}      & $126.74$ & $5.14$ & --- \\
    \emph{One-period policy}& $444.19$ & $322.59$ & --- \\
    \bottomrule
  \end{tabular}
\end{table}

The exact DP policy lowers total cost by $5.14$ hours, or $4.05\%$, relative to the best static allocation, which splits recommendations uniformly at both origins. This gap shows that period-varying allocations pay off even in this small instance. The one-period policy sends all departures to the shorter defaults. It minimizes the current arc-entry cost but ignores the downstream cost of the shared low-capacity arc, and its total travel cost is more than three times the exact minimum. This failure illustrates that a rule that looks best period by period can congest the bottleneck that all later travelers must cross.

Four of five PPO actors reach the minimum cost over the discrete route-share choices, first doing so at validation iterations $30$, $40$, $190$, and $220$, respectively. The selected checkpoints have a mean cost of $121.75$ hours (sample standard deviation $0.33$). The fifth actor's validation cost stabilizes at $122.34$, which exceeds the exact value by $0.61\%$ and improves on the cost of the best static allocation by $3.47\%$. All five actors outperform the best static allocation. Thus, the learning-based approximation recovers the exact optimum in most training runs and remains close to it otherwise, which supports its use on instances where exact computation is out of reach.

\subsection{E2: EMA Bottleneck Performance and Tail Risk}
\label{subsec:exp-e2}

E2 moves to the EMA bottleneck instance, where exact computation is out of reach, and asks whether full-horizon planning pays off relative to simpler rules. E2 uses a Gaussian prior for $\boldsymbol\beta^C$ with means $(0.00,-0.60,0.50,0.50,-0.20)$ and standard deviations $(0.90,0.90,0.70,0.70,0.50)$, in the order of~\eqref{eq:compliance-feature-vector}. We refer to this specification as the biased prior because its means are smaller in magnitude than the true coefficients, and thus it understates travelers' responsiveness to every feature. We evaluate three independently trained versions of the full-horizon PPO policy, each with online belief and recommendation performance updates, together with the nine benchmark policies on $1000$ out-of-sample scenarios with common random numbers. The three PPO actors are training replicates of the same method, whereas the nine benchmarks are alternative routing policies. Online Supplement Section~\ref{app:ema-implementation} details the network construction, stochastic simulator, and PPO training.

Table~\ref{tab:exp-e2-ema} reports, for each routing policy, the mean total cost (column ``Mean''), two cost percentiles (columns ``$\widehat q_{0.90}$'' and ``$\widehat q_{0.99}$''), the upper-$5\%$ average cost (column ``$\widehat{\CVaR}_{0.95}$''), and the paired mean cost difference obtained by subtracting PPO cost from the benchmark-policy cost (column ``Cost minus PPO''). The PPO row averages the three actor-specific estimates of each metric. For each benchmark policy, we compute a $95\%$ Student-$t$ interval from the three training-replicate paired mean cost differences.

\begin{table}[htbp!]
  \centering
  \caption{E2 comparison of PPO with nine benchmark routing policies on the EMA bottleneck instance. 
  }
  \label{tab:exp-e2-ema}
  \scriptsize
  \begingroup
  \setlength{\tabcolsep}{2.7pt}
  \begin{tabular}{lrrrrr}
    \toprule
    Policy & Mean & $\widehat q_{0.90}$ & $\widehat q_{0.99}$ & $\widehat{\CVaR}_{0.95}$ & Cost minus PPO \\
    \midrule
    \emph{PPO full horizon} & $3317.5$ & $3474.1$ & $3602.0$ & $3571.3$ & --- \\
    \emph{Oracle-compliance myopic} & $3504.5$ & $3666.0$ & $3793.5$ & $3770.1$ & $187.0\ [180.6,193.3]$ \\
    \emph{Oracle-compliance-and-road myopic} & $3505.4$ & $3661.6$ & $3787.9$ & $3762.4$ & $187.9\ [181.6,194.2]$ \\
    \emph{Updating myopic} & $3529.0$ & $3692.2$ & $3823.0$ & $3802.7$ & $211.5\ [205.2,217.8]$ \\
    \emph{Static myopic} & $3542.2$ & $3698.6$ & $3850.0$ & $3806.8$ & $224.7\ [218.3,231.0]$ \\
    \emph{Frozen-information myopic} & $3550.0$ & $3715.0$ & $3850.7$ & $3828.1$ & $232.5\ [226.1,238.8]$ \\
    \emph{Predicted-time softmax} & $6901.8$ & $11437.9$ & $13163.7$ & $12794.8$ & $3584.2\ [3577.9,3590.6]$ \\
    \emph{Uniform} & $7980.7$ & $12866.3$ & $15852.2$ & $15123.2$ & $4663.2\ [4656.9,4669.5]$ \\
    \emph{Random Dirichlet} & $8259.5$ & $15097.7$ & $19791.9$ & $18658.0$ & $4942.0\ [4935.6,4948.3]$ \\
    \emph{Default only} & $27019.5$ & $61166.9$ & $65934.5$ & $64929.8$ & $23702.0\ [23695.7,23708.4]$ \\
    \bottomrule
  \end{tabular}
  \endgroup
\end{table}

Three findings stand out. First, PPO has the lowest mean cost and the lowest value of every reported upper-tail statistic among the policies in Table~\ref{tab:exp-e2-ema}. Its three actor-specific mean costs range from $3315.0$ to $3320.1$ hours. Among the myopic benchmarks, \emph{oracle-compliance myopic} has the lowest observed mean cost. Relative to this benchmark, PPO lowers mean cost by $187.0$ hours, or $5.3\%$ of the benchmark policy's mean cost. The $95\%$ interval across training replicates is $[180.6,193.3]$ hours. Relative to the non-oracle \emph{updating myopic} policy, PPO lowers mean cost by $211.5$ hours, or $6.0\%$, with an interval across training replicates of $[205.2,217.8]$ hours. These paired comparisons support PPO's mean cost advantage over both benchmarks across the three training seeds and the common evaluation scenarios.

Second, in this instance, foresight matters more than information. Giving the myopic planner the true compliance vector lowers its mean cost by only $24.5$ hours relative to \emph{updating myopic}, and adding the true road state changes little further. Both oracle benchmarks use coordinate descent to minimize the same full-route proxy for current departures, and thus additional information can change their allocations without improving finite-horizon cost. Their paired mean cost difference is $0.9$ hours, with a $95\%$ confidence interval of $[-1.2,3.1]$ across the common evaluation scenarios. In contrast, PPO does not know the true coefficients and still outperforms both oracle benchmarks by a wide margin. For a planner, this comparison suggests that a forward-looking policy yields larger gains than better behavioral information fed into a myopic rule.

Third, poorly designed guidance is costly, especially in the tails. PPO lowers $\widehat{\CVaR}_{0.95}$ by $5.3\%$ relative to \emph{oracle-compliance myopic}. Relative to \emph{default only}, which reproduces travelers' habitual routes, PPO lowers mean cost by $87.7\%$ and $\widehat{\CVaR}_{0.95}$ by $94.5\%$. The $99$th-percentile costs of \emph{predicted-time softmax}, \emph{uniform}, and \emph{random Dirichlet} range from $13163.7$ to $19791.9$ hours, compared with $3602.0$ hours for PPO. The wide gaps between the mean and upper-percentile costs of these simple rules indicate that guidance that ignores network interactions exposes the operation to large losses in unfavorable scenarios.

Where does PPO's advantage come from? PPO's mean cost consists of $1991.4$ hours of arc-entry travel cost and $1326.1$ hours of terminal penalty, with $203.2$ travelers remaining in the network on average. For \emph{oracle-compliance myopic}, the corresponding values are $2046.8$ hours, $1457.7$ hours, and $229.5$ travelers. The terminal penalty difference accounts for $131.5$ hours of the $187.0$-hour total difference, or $70.4\%$, and the arc-entry travel cost difference accounts for the remaining $55.4$ hours. Because paired policies have identical terminal queue counts, the terminal penalty difference reflects their different numbers of unfinished in-network travelers. In other words, PPO's main advantage is that it brings more travelers to safety within the horizon. Averaging planned route shares across periods, evaluation scenarios, OD groups, and actors, PPO assigns $9.4\%$ to the default route, $73.6\%$ to the first bypass, and $17.0\%$ to the second alternative. \emph{Oracle-compliance myopic} assigns approximately $0.0\%$, $59.9\%$, and $40.1\%$, respectively. PPO thus keeps a modest share on the default route and concentrates most travelers on the first bypass, whereas the myopic benchmark abandons the default route and sends a much larger share to the longer second alternative. One plausible interpretation is that the myopic rule reacts only to the congestion that current departures will create, whereas PPO also accounts for how today's allocation affects travel times and unfinished traffic later in the horizon.

\begin{figure}[htbp!]
  \centering
  \includegraphics[width=0.75\linewidth]{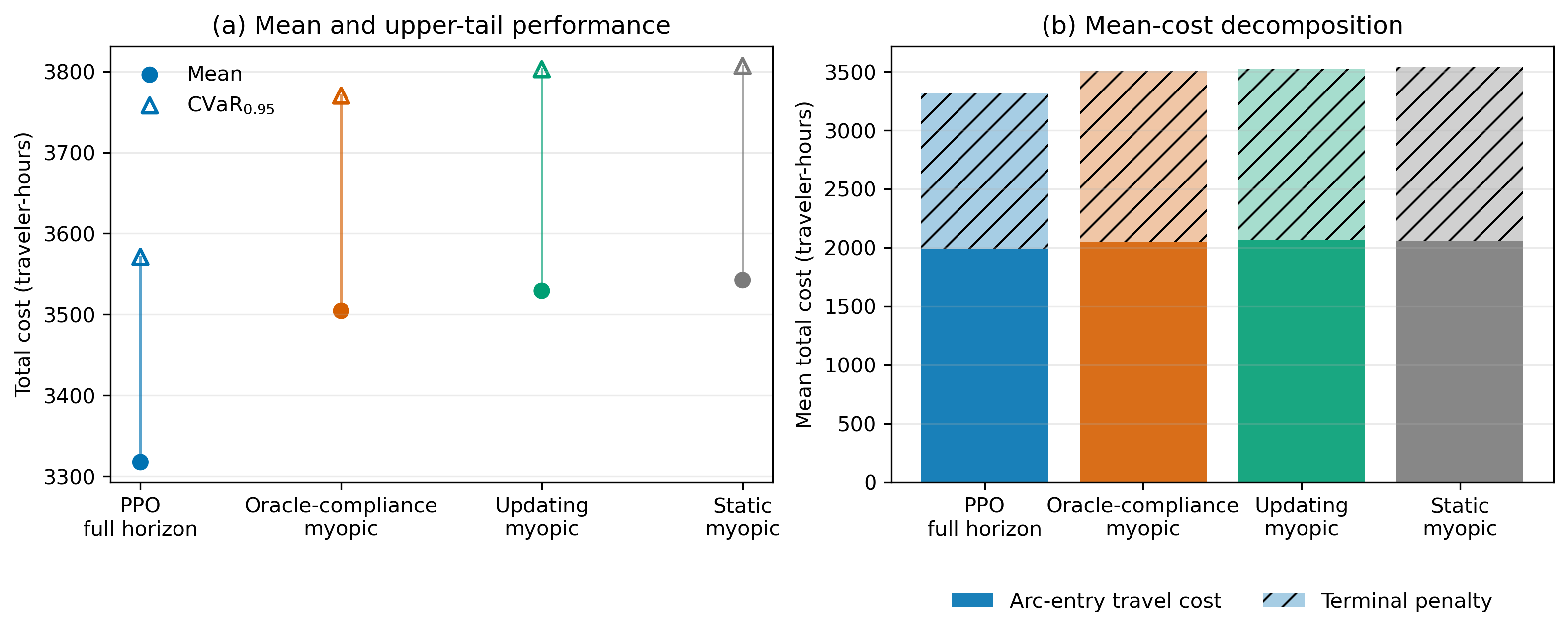}
  \caption{E2 total cost comparison and decomposition for PPO and selected route-aware benchmark policies. Panel~(a) reports mean total cost and $\widehat{\CVaR}_{0.95}$. Panel~(b) decomposes mean total cost into arc-entry travel cost and terminal penalty. 
  }
  \label{fig:exp-e2-cost-profile}
\end{figure}

Figure~\ref{fig:exp-e2-cost-profile} confirms that PPO also incurs less arc-entry cost than every route-aware benchmark policy, and thus the terminal penalty alone does not explain the mean cost ordering. Exogenous departures and arrivals leave every policy with $15.5$ travelers waiting at origins on average, and the terminal penalty differences arise from travelers who remain in the network. 

\subsection{E3: Joint Value of Online Adaptation}
\label{subsec:exp-e3}

E2 shows that full-horizon planning pays off, but it does not isolate the contribution of learning. E3 asks how much the online updates themselves are worth. It retains the EMA bottleneck instance and compares three information settings with separately trained PPO policies. \emph{Frozen} fixes the initial compliance belief, road belief, and recommendation performance throughout the horizon. \emph{Full adaptation} applies the compliance update~\eqref{eq:bayes-compliance-update} through particles~\eqref{eq:particle-weight-update}, the road updates~\eqref{eq:road-filter-update}--\eqref{eq:road-filter-predict}, and the performance update~\eqref{eq:performance-ema} after each period. \emph{Oracle information} receives the true compliance vector and current road state. The \emph{frozen} policy cost minus the \emph{full adaptation} policy cost measures the joint value of compliance updating, road filtering, and performance feedback for these separately trained controllers. The \emph{full adaptation} cost minus the \emph{oracle information} cost indicates how much value online learning forgoes relative to perfect information.

We retain the true compliance vector $\boldsymbol\beta^C=(-0.20,-2.00,1.10,1.00,-0.70)^\top$ and compare three initial Gaussian priors with independent components. We reuse the biased prior from Section~\ref{subsec:exp-e2}. The accurate prior centers near the true coefficients with small standard deviations, and the diffuse prior is weakly informative with large standard deviations. Their means and componentwise standard deviations, in the feature order of~\eqref{eq:compliance-feature-vector}, are
\begin{align*}
\text{Accurate mean:}\quad
&(-0.15,-1.90,1.05,1.00,-0.65),\\
\text{Accurate standard deviation:}\quad
&(0.30,0.35,0.25,0.25,0.25),\\[1mm]
\text{Diffuse mean:}\quad
&(0.00,-0.20,0.20,0.20,0.00),\\
\text{Diffuse standard deviation:}\quad
&(1.40,1.40,1.20,1.20,0.90).
\end{align*}
For each prior and information setting, we use three PPO actors from independent training seeds and a distinct set of $300$ out-of-sample scenarios for each prior and training replicate. Within each set, we evaluate all three information settings on the same scenarios. Online Supplement Section~\ref{app:ema-implementation} specifies the validation and checkpoint selection protocol. The biased-prior \emph{full adaptation} setting uses the same three actors as E2. Because E2 evaluates all three actors on one common set of $1000$ scenarios, these actors have different mean costs in Tables~\ref{tab:exp-e2-ema} and~\ref{tab:exp-e3-adaptation}.

For each prior and training replicate, we average the scenario-level cost differences for two comparisons. We subtract the \emph{full adaptation} cost from the \emph{frozen} cost and the \emph{oracle information} cost from the \emph{full adaptation} cost. Positive differences favor \emph{full adaptation} in the first comparison and \emph{oracle information} in the second. Table~\ref{tab:exp-e3-adaptation} reports mean costs across the three actors in each information setting and $95\%$ Student-$t$ intervals based on the three replicate-level paired differences.

\begin{table}[htbp!]
  \centering
 \caption{E3 comparison of joint online adaptation on the EMA bottleneck instance. We average policy costs across three independently trained actors. Within each prior and training replicate, the three information settings share $300$ out-of-sample evaluation scenarios. Brackets give $95\%$ Student-$t$ intervals for the three training-replicate paired mean differences. 
 }
  \label{tab:exp-e3-adaptation}
  \scriptsize
  \begingroup
  \setlength{\tabcolsep}{3pt}
  \begin{tabular}{lrrrcc}
    \toprule
    Initial prior & \emph{Frozen} & \emph{Full adaptation} & \emph{Oracle information} & Frozen $-$ full & Full $-$ oracle \\
    \midrule
    Accurate & $3412.6$ & $3316.8$ & $3304.0$ & $95.8\ [52.6,139.1]$ & $12.7\ [5.3,20.1]$ \\
    Biased  & $3431.0$ & $3326.4$ & $3315.5$ & $104.5\ [99.9,109.2]$ & $10.9\ [-7.3,29.1]$ \\
    Diffuse & $3427.2$ & $3315.3$ & $3314.3$ & $111.9\ [104.5,119.4]$ & $1.0\ [-27.4,29.4]$ \\
    \bottomrule
  \end{tabular}
  \endgroup
\end{table}

The results deliver two messages. First, learning pays off under every prior that we test. \emph{Full adaptation} lowers mean cost by $2.8\%$, $3.0\%$, and $3.3\%$ relative to \emph{frozen} costs under accurate, biased, and diffuse priors, respectively. All nine within-prior replicate comparisons favor adaptation, and all three intervals exclude zero (Online Supplement Figure~\ref{fig:exp-e3-paired-effects} shows the replicate-level differences). The point estimates increase as the prior moves from accurate to diffuse, although this ordering does not separately identify the effects of prior bias and dispersion. A planner with a poor initial view of traveler behavior thus has at least as much to gain from updating as one with an accurate view. Second, online learning recovers most of the value of perfect information. Mean \emph{full adaptation} cost exceeds \emph{oracle information} cost by $12.7$, $10.9$, and $1.0$ hours under the same priors, respectively, which is small relative to the gains from adaptation. Only the accurate-prior interval excludes zero, and the other intervals establish neither a positive difference nor equality. These separately trained oracle policies do not establish lower bounds on the optimal finite-horizon cost.

Terminal penalties on unfinished in-network travelers account for most of the frozen-minus-full cost differences. Frozen policies leave $16.9$, $19.2$, and $21.2$ additional in-network travelers under accurate, biased, and diffuse priors, respectively. Arc-entry cost differences account for the remaining $11.2$, $8.4$, and $5.8$ hours. Because common exogenous demand gives paired policies identical terminal queues, this decomposition identifies unfinished traffic as the main source of the cost difference, although it does not isolate the contributions of individual belief and performance updates. As in E2, the operational benefit appears mainly as more travelers reaching safety within the horizon.

Because terminal penalties drive most of the cost differences in E2 and E3, we also recompute the saved out-of-sample costs with penalties of $0$, $2.5$, $5$, and $10$ hours per traveler still en route (Online Supplement Section~\ref{app:terminal-sensitivity}). Every point estimate favors PPO or full adaptation at all four penalty values, and the advantages grow with the penalty. Even with no penalty on unfinished traffic, PPO retains a lower mean arc-entry travel cost than \emph{oracle-compliance myopic} in E2, with an interval that excludes zero. A planner who places more weight on clearing the network therefore gains more from dynamic, adaptive recommendations, and a planner who cares only about travel time still benefits from full-horizon planning.

\section{Conclusion}\label{sec:conclusion}

In this paper, we studied dynamic route recommendation under uncertain traveler compliance and road conditions. Our belief-state dynamic program links route recommendations, traffic propagation, compliance learning, road state filtering, and recommendation performance feedback, and we showed that a compact belief state suffices for optimal control. A particle approximation and a PPO actor with a Dirichlet policy provide a computational solution for continuous route shares.
The experiments support dynamic recommendations and online updating, and they suggest three managerial insights. First, foresight matters. In the deterministic instance, PPO attains the minimum cost over the discrete route-share choices in four of five training runs. On the EMA bottleneck instance, full-horizon PPO reduces mean total cost and the upper-$5\%$ average cost by about $5.3\%$ relative to the myopic benchmark with the lowest mean cost, even though that benchmark knows the true compliance coefficients. Second, the main benefit of planning ahead is clearing the network. Fewer travelers remaining en route explain most of the mean cost reduction, and PPO also incurs lower arc-entry travel costs. Third, learning compensates for imperfect knowledge of travelers. Updating beliefs and recommendation performance lowers mean total cost under accurate, biased, and diffuse initial priors, and the remaining gap to policies with perfect information is small. 

For future research, one can test sparse telemetry, different true compliance coefficients, and destination choices with binding receiving capacities. 
One can also characterize when the planner should deliberately issue informative recommendations to accelerate learning, which is the dual-control trade-off that Remark~\ref{remark:1} highlights. 


\bibliographystyle{informs2014} 
\bibliography{literature} 

\clearpage
\numberwithin{equation}{section}
\numberwithin{table}{section}
\numberwithin{figure}{section}

\begin{appendices}

\section{Technical Derivations and Robustness Extensions}
\label{app:technical-results}

We present two extensions of the response model, derive the road filter, and provide details proofs for Propositions~\ref{prop:belief-factorization} and~\ref{prop:belief-state-sufficiency}. Throughout, $b_t^C$ is a probability measure on $\mathbb R^5$, possibly without a Lebesgue density. Conditional identities hold almost surely, and posterior formulas require positive predictive likelihood.

\subsection{Multiple Destination-Consistent Rejection Routes}
\label{app:multiple-default-routes}

For each informative route $p\in\mathcal P^{\mathrm{inf}}_{o,t}(S_t)$, let $\mathcal Q^0_{op,t}(S_t)\subseteq\mathcal P_{o,d(p),t}(S_t)\setminus\{p\}$ be the nonempty set of alternatives that a traveler may select after rejecting $p$. The traveler selects $q\in\mathcal Q^0_{op,t}(S_t)$ with probability $\lambda^0_{opq,t}$, a known nonnegative measurable function of $S_t$ that sums to one over $q$. The nominal default route $p^0_{op,t}$ remains the reference for compliance features and recommendation performance, with deterministic agreement when the recommendation coincides with it.

Let $Z_{opq,t}$ count travelers receiving recommendation $p$ and selecting alternative $q$, and collect these counts in $\mathbf Z_t$. Conditional independence of traveler responses given pre-decision information and the compliance coefficient vector gives
\begin{equation}\label{eq:app-multiple-default-response}
    \begin{aligned}
    \left(
        Y_{op,t},(Z_{opq,t})_{q\in\mathcal Q^0_{op,t}(S_t)}
    \right)
    \mid\mathscr F_t,\mathbf r_t,\mathbf N_t,\boldsymbol\beta^C
    \sim
    \operatorname{Multinomial}\!\left(
        N_{op,t};
        \pi^C_{op,t},
        \bigl((1-\pi^C_{op,t})\lambda^0_{opq,t}\bigr)_q
    \right).
    \end{aligned}
\end{equation}
The independence assumptions in Section~\ref{subsec:dp-compliance-learning} still apply. Because $\mathcal Q^0_{op,t}(S_t)$ excludes $p$, a traveler assigned to route $p$ who selects $p$ must have accepted the recommendation. Thus, $Y_{op,t}$ follows~\eqref{eq:compliance-binomial}. Given $(S_t,N_{op,t},Y_{op,t})$, the alternative counts follow a multinomial distribution with total $N_{op,t}-Y_{op,t}$ and probability $\lambda^0_{opq,t}$ for each $q \in \mathcal Q^0_{op,t}(S_t)$.

The planner must observe each traveler's selected route by the end of the departure period. We replace the default-route contribution $(N_{op,t}-Y_{op,t})\mathbf 1\{e_{p^0_{op,t},1}=(i,j)\}$ in~\eqref{eq:dp-arc-flow} by $\sum_{q\in\mathcal Q^0_{op,t}(S_t)} Z_{opq,t}\mathbf 1\{e_{q,1}=(i,j)\}$ and propagate packets on their selected complete routes. We augment $O^C_{t+1}$ by $\mathbf Z_t$ and include these counts in $\mathbf c_t^\Theta$ and $\Xi_{t+1}$ in Equations~\eqref{eq:road-filter-context} and~\eqref{eq:random-information}. For known $\lambda^0$, Equation~\eqref{eq:bayes-compliance-update} remains valid because the likelihood factor for $\mathbf Z_t$ conditional on $(S_t,\mathbf r_t,\mathbf N_t,\mathbf Y_t)$ is independent of $\boldsymbol\beta^C$. Equations~\eqref{eq:road-filter-update}--\eqref{eq:road-filter-predict} and~\eqref{eq:bellman-main} keep their forms under these changes.
If $\lambda^0$ is unknown, learning it requires a joint posterior with $\boldsymbol\beta^C$ and a likelihood that uses the alternative counts. If a rejecting traveler could also select $p$, observed agreement would no longer identify acceptance, and~\eqref{eq:compliance-likelihood} would require a different likelihood. Setting $\mathcal Q^0_{op,t}(S_t)=\{p^0_{op,t}\}$ recovers the baseline.

\subsection{Beta--Binomial Robustness Extension}
\label{app:beta-binomial}

To allow positive dependence among responses within an informative group--route--period cell, we replace the baseline conditional independence assumption with a cell-level latent agreement probability and draw it anew each period with precision $\phi_{\mathrm{BB}}>0$:
\begin{align}
\label{eq:beta-binomial-latent-probability}
    \widetilde\pi^C_{op,t}
    \mid\mathscr F_t,\boldsymbol\beta^C,\phi_{\mathrm{BB}}
    &\sim
    \operatorname{Beta}\!\left(
        \phi_{\mathrm{BB}}\pi^C_{op,t},
        \phi_{\mathrm{BB}}(1-\pi^C_{op,t})
    \right),\\
\label{eq:beta-binomial-conditional-count}
    Y_{op,t}\mid
        \mathscr F_t,\mathbf N_t,\widetilde\pi^C_{op,t}
    &\sim
    \operatorname{Binomial}
        (N_{op,t},\widetilde\pi^C_{op,t}).
\end{align}
Given $(\mathscr F_t,\boldsymbol\beta^C,\phi_{\mathrm{BB}})$, we draw the latent probabilities independently across cells. Further conditioning on the action, assignment counts, road state, or new demand counts leaves their law unchanged. Given the history, assignments, and all cell probabilities, agreement counts independently follow~\eqref{eq:beta-binomial-conditional-count}, with no further dependence on $\boldsymbol\beta^C$, $\phi_{\mathrm{BB}}$, the action, road state, or new demand counts. These assumptions preserve conditional independence between behavioral parameters and road state.

Marginalizing over $\widetilde\pi^C_{op,t}$ gives
\begin{equation}
\label{eq:beta-binomial-count}
    \begin{aligned}
    Y_{op,t}\mid
        \mathscr F_t,\mathbf N_t,\boldsymbol\beta^C,\phi_{\mathrm{BB}}
    \sim
    \operatorname{BetaBinomial}\!\left(
        N_{op,t},
        \phi_{\mathrm{BB}}\pi^C_{op,t},
        \phi_{\mathrm{BB}}(1-\pi^C_{op,t})
    \right).
    \end{aligned}
\end{equation}
Conditional on the variables in~\eqref{eq:beta-binomial-count}, the mean remains $N_{op,t}\pi^C_{op,t}$, the variance inflates the binomial variance by the factor $1+(N_{op,t}-1)/(\phi_{\mathrm{BB}}+1)$, and two distinct travelers' agreement indicators in the same cell have correlation $1/(\phi_{\mathrm{BB}}+1)$.

Writing $\mathrm{B}(\cdot,\cdot)$ for the beta function and setting $a^{\mathrm{BB}}_{op,t}=\phi_{\mathrm{BB}}\pi^C_{op,t}$ and $b^{\mathrm{BB}}_{op,t}=\phi_{\mathrm{BB}}(1-\pi^C_{op,t})$, we obtain the likelihood
\begin{equation}
\label{eq:beta-binomial-likelihood}
    \begin{aligned}
    L_t^{C,\mathrm{BB}}(\boldsymbol\beta^C,\phi_{\mathrm{BB}})
     =
    \prod_{o\in\mathcal O_t^+(S_t)}
    \prod_{p\in\mathcal P^{\mathrm{inf}}_{o,t}(S_t)}
    \left[
    \binom{N_{op,t}}{Y_{op,t}}
    \frac{
        \mathrm{B}\left(Y_{op,t}+a^{\mathrm{BB}}_{op,t},
            N_{op,t}-Y_{op,t}+b^{\mathrm{BB}}_{op,t}\right)
    }{
        \mathrm{B}\left(a^{\mathrm{BB}}_{op,t},b^{\mathrm{BB}}_{op,t}\right)
    }
    \right].
    \end{aligned}
\end{equation}
Smaller $\phi_{\mathrm{BB}}$ gives stronger dependence and greater overdispersion. As $\phi_{\mathrm{BB}}\to\infty$, Equation~\eqref{eq:beta-binomial-count} converges to the binomial law in~\eqref{eq:compliance-binomial}. Cells with $p=p^0_{op,t}$ remain deterministic under~\eqref{eq:coincident-route-agreement} and do not enter~\eqref{eq:beta-binomial-likelihood}.

For fixed $\phi_{\mathrm{BB}}$, replacing $L_t^C$ with $L_t^{C,\mathrm{BB}}$ in~\eqref{eq:bayes-compliance-update} gives the exact update. Learning $\phi_{\mathrm{BB}}$ requires a joint posterior on $(\boldsymbol\beta^C,\phi_{\mathrm{BB}})$, and Proposition~\ref{prop:belief-factorization} continues to hold for this enlarged parameter if it is independent of $\Theta_0$ given $\mathscr F_0$. Because the latent cell probabilities do not persist across periods, we integrate them out before forming the next-period state.

\subsection{Proof of Conditional Belief Factorization}
\label{app:belief-factorization-proof}
\label{app:compliance-posterior}

\proof{Proof of Proposition~\ref{prop:belief-factorization}.}
We argue by induction. Equation~\eqref{eq:initial-belief-factorization} gives the base case $t=0$. Suppose~\eqref{eq:belief-factorization} holds at epoch $t$. The demand and assignment assumptions in Section~\ref{subsec:dp-demand-counts} give, with $\mathcal L$ denoting a conditional law,
\begin{equation}
\label{eq:app-demand-exogeneity}
    \mathcal L\!\left(
        \mathbf D_t^{\mathrm{act}},\mathbf D_{t+1}^{\mathrm{new}}
        \mid\mathscr F_t,\mathbf r_t,\boldsymbol\beta^C,\Theta_t
    \right)
    =
    K_t^D(\cdot\mid\mathbf q_t)
\end{equation}
and
\begin{equation}
\label{eq:app-assignment-exogeneity}
    \begin{aligned}
    \mathcal L\!\left(
        \mathbf N_t
        \mid\mathscr F_t,\mathbf r_t,
        \mathbf D_t^{\mathrm{act}},\mathbf D_{t+1}^{\mathrm{new}},
        \boldsymbol\beta^C,\Theta_t
    \right)
     =
    \bigotimes_{o\in\mathcal O_t^+(S_t)}
    \operatorname{Multinomial}
        (D_{o,t}^{\mathrm{act}},\mathbf r_{o,t}).
    \end{aligned}
\end{equation}
Assignments are independent across active groups. When no group is active, the empty count vector has probability one. Conditional on $\mathscr F_t$, policy randomization is independent of $\boldsymbol\beta^C$ and $\Theta_t$. Conditioning on the action therefore leaves their beliefs unchanged. We fix the action and the realized demand and assignments. The features $\boldsymbol\varphi_{op,t}$ are $\mathscr F_t$-measurable and fixed across candidate coefficients. Equations~\eqref{eq:compliance-binomial}--\eqref{eq:coincident-route-agreement} give the likelihood $L_t^C(\boldsymbol\beta^C)$ for $\mathbf Y_t$ on its support. Further conditioning on $\Theta_t$ and $\mathbf D_{t+1}^{\mathrm{new}}$ leaves this likelihood unchanged.

Travelers report selected routes by their departure period's end (Section~\ref{subsec:dp-compliance-learning}), and thus observing agreement does not require trip completion. Given $(\Theta_t,\mathbf c_t^\Theta)$, the road observation law~\eqref{eq:road-observation-kernel} has no further dependence on $\boldsymbol\beta^C$, $\mathbf D_{t+1}^{\mathrm{new}}$, or $\mathscr F_t$. Thus, apart from parameter-independent factors, the observation likelihood separates into $L_t^C(\boldsymbol\beta^C)$ and $L_t^\Theta(O^\Theta_{t+1}\mid\vartheta,\mathbf c_t^\Theta)$. Bayes' rule gives, for measurable $B\subseteq\mathbb R^5$,
\begin{equation}
\label{eq:app-period-posterior-factorization}
    \begin{aligned}
    \Pr\left(
        \boldsymbol\beta^C\in B,\Theta_t=\vartheta
        \mid\mathscr F_{t+1}
    \right)
     =
    \frac{
        \displaystyle\int_B L_t^C(\boldsymbol\beta^C)
            b_t^C(d\boldsymbol\beta^C)
    }{
        \displaystyle\int_{\mathbb R^5}
            L_t^C(\boldsymbol\beta^C)b_t^C(d\boldsymbol\beta^C)
    }
    \frac{
        L_t^\Theta(O^\Theta_{t+1}\mid\vartheta,\mathbf c_t^\Theta)
        b_t^\Theta(\vartheta)
    }{
        \displaystyle\sum_{u\in\mathcal S^\Theta}
        L_t^\Theta(O^\Theta_{t+1}\mid u,\mathbf c_t^\Theta)
        b_t^\Theta(u)
    }.
    \end{aligned}
\end{equation}
With observed context fixed, the product prior and the separated likelihoods factor the normalizing constant. The demand and assignment factors cancel. Cells with $p=p^0_{op,t}$ contribute a factor of one on their support. Summing over $\vartheta$ gives~\eqref{eq:bayes-compliance-update}. If all informative counts vanish, $L_t^C\equiv1$ and $b_{t+1}^C=b_t^C$. Hence the joint posterior is $b_{t+1}^C(B)\widetilde b_t^\Theta(\vartheta)$, where $\widetilde b_t^\Theta(\vartheta)=\Pr(\Theta_t=\vartheta\mid\mathscr F_{t+1})$. The exogenous road transition~\eqref{eq:road-transition-kernel} then gives
\begin{align}
\label{eq:app-next-period-factorization}
    \Pr\!\left(
        \boldsymbol\beta^C\in B,\Theta_{t+1}=\vartheta'
        \mid\mathscr F_{t+1}
    \right)\notag
    &=
    \sum_{\vartheta\in\mathcal S^\Theta}
        \Pr(\Theta_{t+1}=\vartheta'\mid\Theta_t=\vartheta)
        \Pr\!\left(
            \boldsymbol\beta^C\in B,\Theta_t=\vartheta
            \mid\mathscr F_{t+1}
        \right)\notag\\
    & =
    b_{t+1}^C(B)
    \sum_{\vartheta\in\mathcal S^\Theta}
        \Pr(\Theta_{t+1}=\vartheta'\mid\Theta_t=\vartheta)
        \widetilde b_t^\Theta(\vartheta)
    =
    b_{t+1}^C(B)b_{t+1}^\Theta(\vartheta').
\end{align}
This establishes~\eqref{eq:belief-factorization} at epoch $t+1$ and completes the induction. \Halmos
\endproof

\subsection{Road-Filter Derivation}
\label{app:road-filter}

Marginalizing~\eqref{eq:app-period-posterior-factorization} over $\boldsymbol\beta^C$ gives~\eqref{eq:road-filter-update}. Equation~\eqref{eq:app-next-period-factorization} gives~\eqref{eq:road-filter-predict}. These equations distinguish the completed-period belief $\widetilde b_t^\Theta$ from the next-decision belief $b_{t+1}^\Theta$. Conditional on $(\mathscr F_t,\mathbf r_t)$, demand, assignment, and agreement counts carry no information about $\Theta_t$ by exogeneity and Proposition~\ref{prop:belief-factorization}. The realized counts enter $\mathbf c_t^\Theta$ because they determine traffic loading.

The likelihood~\eqref{eq:road-observation-likelihood} must use all retained road information, and thus observed packet transitions, entry durations, and costs belong in $O^\Theta_{t+1}$ or its sufficient statistic. Because we define it relative to a dominating measure, it covers discrete, continuous, mixed, and deterministic measurements, and it factors into component likelihoods only if the components are conditionally independent given $(\Theta_t,\mathbf c_t^\Theta)$.

If factorization fails, counts may convey road information, and the separate updates~\eqref{eq:bayes-compliance-update} and~\eqref{eq:road-filter-update}--\eqref{eq:road-filter-predict} generally lose sufficiency. A joint posterior over $(\boldsymbol\beta^C,\Theta_t)$ suffices only if this pair contains all latent variables governing future transitions. Otherwise, the state must also track the remaining persistent latent variables.

\subsection{Proof of Belief-State Sufficiency}
\label{app:belief-state-proof}

\proof{Proof of Proposition~\ref{prop:belief-state-sufficiency}.}
By Proposition~\ref{prop:belief-factorization}, the conditional law of $(\boldsymbol\beta^C,\Theta_t)$ given $\mathscr F_t$ is $b_t^C\,b_t^\Theta$. Given an admissible action, demand and assignments depend on history only through $(S_t,\mathbf r_t)$ by~\eqref{eq:app-demand-exogeneity}--\eqref{eq:app-assignment-exogeneity}. Agreement counts additionally depend on $\boldsymbol\beta^C$ and on the state-determined features through~\eqref{eq:compliance-binomial}--\eqref{eq:coincident-route-agreement}. Traffic loading and road observations additionally use $\Theta_t$ and realized counts. Integrating these laws against $b_t^C(d\boldsymbol\beta^C)b_t^\Theta(\vartheta)$ gives~\eqref{eq:predictive-distribution}. The exogenous transition~\eqref{eq:road-transition-kernel} determines the next road belief through~\eqref{eq:road-filter-predict}.

Because the packet records retain each packet's route, position, and residual duration, they determine which inherited packets enter which arcs in period $t$. Given $(S_t,\mathbf r_t,\Xi_{t+1},\mathbf O_{t+1})$, the deterministic function $f_t$ in Section~\ref{subsec:dp-transition} determines $S_{t+1}$, and Equation~\eqref{eq:dp-stage-cost} determines $g_t$. For every bounded measurable function $\chi$,
\begin{equation}\label{eq:app-markov-prediction}
    \begin{aligned}
    \mathbb E\!\left[
        \chi(S_{t+1},g_t)
        \mid\mathscr F_t,\mathbf r_t
    \right]
     =
    \int
    \chi\!\left(
        f_t(S_t,\mathbf r_t,\zeta,\boldsymbol\xi),
        g_t(S_t,\mathbf r_t,\zeta)
    \right)
    P_t(d\zeta,d\boldsymbol\xi\mid S_t,\mathbf r_t).
    \end{aligned}
\end{equation}
The right-hand side of~\eqref{eq:app-markov-prediction} depends on history only through $(S_t,\mathbf r_t)$, proving the controlled Markov property. Because $G_T$ is measurable, backward induction under the regularity and attainment assumptions in Section~\ref{subsec:dp-transition} yields Equation~\eqref{eq:bellman-main} and an optimal measurable Markov selector. The proof establishes sufficiency for the exact state $S_t$, which includes complete packet records and exact posterior measures. Establishing sufficiency for the computational state $S_t^{(M)}$ would require a separate argument. \Halmos
\endproof

\section{Additional Experimental and Implementation Details}
\label{app:additional-experiments}

\subsection{PPO Training Updates}
\label{app:ppo-training}

This section states the training updates that Section~\ref{subsec:ppo-reward-reformulation} summarizes. The subscript $\mathrm{fixed}$ identifies the parameter values $\boldsymbol\omega_{\mathrm{fixed}}$ and $\boldsymbol\nu_{\mathrm{fixed}}$ at the start of a roll-out batch. We use the actor at these values to generate the batch and the critic at these values to compute its advantages and critic targets. Before optimization, we record the log densities of the sampled shares and the critic values and hold them fixed while updating the candidate parameters.
We compute temporal-difference residuals with the fixed critic $V_{\boldsymbol\nu_{\mathrm{fixed}}}$:
\begin{equation}
\label{eq:td-residual}
    \delta_t
    =
    R_{t+1}
    +
    \gamma V_{\boldsymbol\nu_{\mathrm{fixed}}}(\widetilde S_{t+1})
    -
    V_{\boldsymbol\nu_{\mathrm{fixed}}}(\widetilde S_t), 
\end{equation}
and compute the generalized advantage estimates as
\begin{equation}
\label{eq:gae}
    \widehat A_t
    =
    \sum_{\ell=0}^{T-t-1}
    (\gamma\lambda)^\ell\delta_{t+\ell},
\end{equation}
where $\lambda\in[0,1]$ weights subsequent residuals. 
Because Equation~\eqref{eq:ppo-reward} includes $G_T$ in the final reward and sets the terminal critic value to zero, it counts the terminal penalty exactly once.

For each recorded actor share vector $\mathbf r_t^{\mathrm{actor}}$ and feature vector $\widetilde S_t$, we compute the policy likelihood ratio
\begin{equation}
\label{eq:ppo-ratio}
    \rho_t(\boldsymbol\omega)
    =
    \frac{
        \pi_{\boldsymbol\omega}(\mathbf r_t^{\mathrm{actor}}\mid\widetilde S_t)
    }{
        \pi_{\boldsymbol\omega_{\mathrm{fixed}}}(\mathbf r_t^{\mathrm{actor}}\mid\widetilde S_t)
    }.
\end{equation}
The numerator uses the candidate actor parameters $\boldsymbol\omega$, and the denominator uses the fixed parameters $\boldsymbol\omega_{\mathrm{fixed}}$ that generated the recorded shares. We evaluate both densities on the same recorded shares, feature vector, and available route sets. For clipping tolerance $\varepsilon_{\mathrm{PPO}}>0$, the PPO surrogate for a single sample is
\begin{align}
\label{eq:ppo-clip-objective}
    \ell_t^{\mathrm{CLIP}}(\boldsymbol\omega)
    =
    \min\Bigl\{
        \rho_t(\boldsymbol\omega)\widehat A_t,
        \operatorname{clip}\!\left(
            \rho_t(\boldsymbol\omega),
            1-\varepsilon_{\mathrm{PPO}},
            1+\varepsilon_{\mathrm{PPO}}
        \right)\widehat A_t
    \Bigr\},
\end{align}
where $\operatorname{clip}(x,a,b)$ truncates $x$ to the interval $[a,b]$. We update the parameters of the shared hidden layers and both output heads by maximizing
\begin{equation}
\label{eq:ppo-full-objective}
    \mathbb E_t\!\left[
        \ell_t^{\mathrm{CLIP}}(\boldsymbol\omega)
        -
        c_V\left(
            V_{\boldsymbol\nu}(\widetilde S_t)-\widehat R_t
        \right)^2
        +
        c_H\mathsf{Ent}\!\left[
            \pi_{\boldsymbol\omega}(\cdot\mid\widetilde S_t)
        \right]
    \right],
\end{equation}
where $\mathbb E_t$ denotes the sampled batch average and $\mathsf{Ent}$ sums the differential entropies of the Dirichlet factors in~\eqref{eq:joint-dirichlet-policy}. 
The coefficients $c_V,c_H\geq0$ weight the critic loss and policy entropy. Before optimization, we compute the critic target $\widehat R_t=\widehat A_t+V_{\boldsymbol\nu_{\mathrm{fixed}}}(\widetilde S_t)$ and hold $\widehat R_t$ and $\widehat A_t$ fixed during batch updates. The clipping in~\eqref{eq:ppo-clip-objective} keeps each update close to the policy that generated the batch, which stabilizes training when simulated costs are noisy. Positive reward rescaling preserves the unregularized control objective up to scale but changes the relative influence of the entropy and critic-loss terms unless we adjust their coefficients. 

\subsection{E1: Exact Validation with Discrete Route Shares}
\label{app:e1-details}

\subsubsection{Instance and objective.}

Table~\ref{tab:app-e1-network} lists each origin's default route, which passes through $(A,D)$, and its detour, which passes through $(B,D)$. All arcs use~\eqref{eq:dp-travel-time} with $\alpha_{\mathrm{BPR}}=0.15$, $\beta_{\mathrm{BPR}}=4$, and unit-length periods. Times use arbitrary units, and costs use traveler-time units and unrounded arc times. Packet propagation uses the rounded durations in~\eqref{eq:packet-duration}.

\begin{table}[htbp!]
  \centering
  \caption{E1 route and arc calibration. Entries in parentheses follow the listed arc sequence.}
  \label{tab:app-e1-network}
  \small
  \begin{tabular}{llrr}
    \toprule
    Route & Arc sequence & Free-flow time & Capacity \\
    \midrule
    $O_1$ default & $(O_1,A)$, $(A,D)$ & $(1.0,2.0)$ & $(8,4)$ \\
    $O_2$ default & $(O_2,A)$, $(A,D)$ & $(1.0,2.0)$ & $(8,4)$ \\
    $O_1$ detour  & $(O_1,B)$, $(B,D)$ & $(1.4,3.0)$ & $(8,9)$ \\
    $O_2$ detour  & $(O_2,B)$, $(B,D)$ & $(1.4,3.0)$ & $(8,9)$ \\
    \bottomrule
  \end{tabular}
\end{table}

Each origin starts with $15$ waiting travelers and no new demand. Departures are $(8,2,5)$ for $O_1$ and $(2,8,5)$ for $O_2$. Largest-remainder apportionment rounds each share-based count down and then assigns the remaining travelers in decreasing order of fractional remainder. Every traveler follows the recommended route, and agreement counts therefore satisfy $Y_{op,t}=N_{op,t}$ deterministically. Identical, equally weighted particles at $\boldsymbol\beta^C=(12,0,0,0,0)^\top$ supply fixed compliance inputs. We hold the road state at its normal value throughout the planning horizon and disable compliance, road-belief, and recommendation performance updates.

Each origin's default-route share $x_{o,t}=r_{op^0,t}$ takes a value in $\{0,0.5,1\}$ and determines its detour share $1-x_{o,t}$. Here $p^0$ denotes that origin's unique default route. The two origins' shares yield nine joint actions $\mathcal A_G=\{0,0.5,1\}^2$. The three-period objective is
\begin{equation}
\label{eq:app-e1-objective}
  Z^\pi_{\mathrm{E1}}
  =\sum_{t=0}^{2}\sum_{(i,j)\in\mathcal E}
    z_{ij,t}\operatorname{tt}_{ij}(z_{ij,t})
    +W_{15}(\mathbf X_3),
\end{equation}
where $\operatorname{tt}_{ij}$ gives normal-state arc times. The term $W_{15}(\mathbf X_3)$ replaces the baseline terminal penalty. It sums arc-entry costs as the network clears $\mathbf X_3$ over periods $3$--$17$ with no new departures. Every policy clears the network within this continuation.

The Bellman recursion over $\mathcal A_G$ minimizes cost over all $9^3=729$ three-period schedules from the fixed initial state. We also evaluate all nine static schedules that repeat one action. The one-period policy minimizes current arc-entry cost, including inherited entries, over the same action set. For this instance, the static schedule with uniform route shares, $(x_{O_1,t},x_{O_2,t})=(0.5,0.5)$ in every period, attains the lowest cost among the nine static schedules. The one-period policy selects $(1,1)$ in every period and therefore coincides with the policy that recommends only the default routes.

Every occupied first arc has a rounded duration of two periods, and thus travelers departing in period $t$ enter the shared downstream arcs in period $t+2$. Total cost therefore separates into the three cohorts' complete route costs, and recovering an exact schedule tests whether PPO can optimize under delayed cost feedback.

\subsubsection{Categorical PPO and training protocol.}

In E1, the actor samples categorical actions from $\mathcal A_G$ during training and selects the most probable action during validation and evaluation. Categorical log probabilities and entropy replace the Dirichlet quantities in~\eqref{eq:joint-dirichlet-policy}, \eqref{eq:ppo-ratio}, and~\eqref{eq:ppo-full-objective}. We retain Section~\ref{sec:ppo-approximation}'s state encoding, critic, advantage estimation, and parameter updates, using E1 dynamics and objective~\eqref{eq:app-e1-objective}.

The actor and critic share two fully connected $128$-unit hidden layers with hyperbolic-tangent activations, and separate output heads produce nine action logits and one value. The inputs follow the E2--E3 encoding in Section~\ref{app:ema-ppo-implementation}, and the belief and performance inputs remain constant because E1 disables these updates.

We train five independent actors for $240$ iterations each on CPU. Each iteration collects $360$ transitions ($120$ complete trajectories) and performs $5$ update epochs with minibatches of $180$. The remaining settings are a learning rate of $3\times10^{-4}$, a PPO clip radius of $0.20$, an entropy coefficient of $0.01$, a value-loss coefficient of $0.50$, $\gamma=\lambda=1$, a reward scale of $0.01$, and a reward clip of $10$.

We validate the argmax policy after iteration $1$, every $10$ iterations, and at iteration $240$, and we treat a cost within $10^{-6}$ of the exact minimum as a recovery. We save a checkpoint whenever the validation cost strictly improves. We evaluate each restored actor $20$ times using actor-specific seeds. Because both the argmax policy and E1 dynamics are deterministic, all $20$ checks for a given actor must yield the same total cost. The five independent training runs measure optimization variability.

\subsubsection{Actor-level recovery and action schedules.}

Four of the five training runs (seeds $10{,}123$, $20{,}123$, $30{,}123$, and $40{,}123$, with selected iterations $40$, $220$, $190$, and $30$) attain the exact minimum of $121.60$ within $10^{-6}$ traveler-time units. Actor $0$ (seed $123$, selected iteration $30$) attains $122.34$ and exceeds the minimum by $0.74$. Across the five actors, the mean cost is $121.75$ traveler-time units and the sample standard deviation is $0.33$ (Section~\ref{subsec:exp-e1}).

Actors $1$--$4$ recover two cost-equivalent exact schedules (Table~\ref{tab:app-e1-actions}). The period-$2$ actions $(0,1)$ and $(1,0)$ have equal cost because both origins release five travelers in period $2$ under identical arc calibrations.

\begin{table}[htbp!]
  \centering
  \caption{E1 selected action schedules. Each pair is $(x_{O_1,t},x_{O_2,t})$, where $x_{o,t}$ is origin $o$'s default-route share.}
  \label{tab:app-e1-actions}
  \small
  \begin{tabular}{lccc}
    \toprule
    Policy & Period 0 & Period 1 & Period 2 \\
    \midrule
    \policy{Exact discrete-action DP} & $(0.5,0)$ & $(0,0.5)$ & $(0,1)$ \\
    \policy{PPO actor 0} & $(0.5,0.5)$ & $(0,0.5)$ & $(1,0)$ \\
    \policy{PPO actor 1} & $(0.5,0)$ & $(0,0.5)$ & $(1,0)$ \\
    \policy{PPO actor 2} & $(0.5,0)$ & $(0,0.5)$ & $(1,0)$ \\
    \policy{PPO actor 3} & $(0.5,0)$ & $(0,0.5)$ & $(0,1)$ \\
    \policy{PPO actor 4} & $(0.5,0)$ & $(0,0.5)$ & $(0,1)$ \\
    \policy{Best static allocation} & $(0.5,0.5)$ & $(0.5,0.5)$ & $(0.5,0.5)$ \\
    \policy{One-period policy} & $(1,1)$ & $(1,1)$ & $(1,1)$ \\
    \bottomrule
  \end{tabular}
\end{table}

Arc-entry costs during the decision periods are $64.99$ for the exact schedule, $65.73$ for actor $0$, $64.68$ for the best static allocation, and $169.82$ for the one-period policy. The corresponding clearance costs are $56.61$, $56.61$, $62.06$, and $274.38$. Relative to the best static allocation, the exact schedule therefore incurs $0.31$ more decision-period cost and $5.45$ less clearance cost, saving $5.14$ overall. Downstream clearance entries determine the ranking.

Actor $0$'s period-$0$ action $(0.5,0.5)$ costs $12.09$, versus $12.49$ for the exact action $(0.5,0)$. Its downstream entries in period $2$ cost $41.15$, versus $40.01$ under the exact schedule. The delayed increase exceeds the immediate saving by $0.74$. This gap equals actor $0$'s total cost difference at every validation from iteration $30$ through iteration $240$.

\subsection{Implementation Details for E2 and E3}
\label{app:ema-implementation}

E2 and E3 share the following simulator and PPO implementation. Sections~\ref{subsec:exp-e2}--\ref{subsec:exp-e3} specify their compliance priors and information settings.

\subsubsection{Network extraction and route catalog.}

The EMA network has $74$ nodes and $258$ directed links. Groups $\mathcal O=\{o_{30,23},o_{30,24},o_{31,23},o_{31,24}\}$ occupy two physical origins, with $v(o_{a,d})=a$ and fixed destinations $\mathcal D_{o_{a,d}}=\{d\}$. We use Yen's algorithm~\citep{yen1971finding} with deterministic tie-breaking and Dijkstra's algorithm~\citep{dijkstra1959note} for its subproblems. We compute each group's three shortest loop-less paths using the fixed nominal arc times $\bar t^0_{ij,t}=t^0_{ij,t}(\mathrm{normal})$. The resulting route catalog spans $13$ links.

Within each group, we index the routes $p^0,p^1,p^2$ in increasing order of free-flow time. The first route $p^0$ is the destination-consistent default $p^0_{op,t}$ for every recommendation $p$ in the group. Because we lack route-specific familiarity data, we set $\mathrm{Fam}_{op^0}=0.95$ and specify the other scores as
\begin{equation}
\label{eq:app-route-familiarity}
    \mathrm{Fam}_{op^j}
    =
    \left[
        0.95-0.16j
        -0.22\frac{
            \operatorname{TT}^{0}_{p^j}-\operatorname{TT}^{0}_{p^0}
        }{\operatorname{TT}^{0}_{p^0}}
        -0.03\max\{|p^j|-|p^0|,0\}
    \right]_{0.35}^{0.85},
    \quad j\in\{1,2\},
\end{equation}
where $\operatorname{TT}^{0}_{p}$ denotes free-flow time and $[x]_a^b=\min\{b,\max\{a,x\}\}$. These scores depend on route rank, proportional excess free-flow time, and excess arc count. Table~\ref{tab:app-ema-route-catalog} reports the resulting scores.

\begin{table}[htbp!]
  \centering
  \caption{Route catalog for E2 and E3. Within each group, $p^0$ is
  the destination-consistent default, $p^1$ the first bypass, and
  $p^2$ the second alternative.}
  \label{tab:app-ema-route-catalog}
  \begingroup
  \footnotesize
  \setlength{\tabcolsep}{4pt}
  \renewcommand{\arraystretch}{1.12}
  \begin{tabularx}{\linewidth}{
    @{}lc>{\raggedright\arraybackslash}Xrrc@{}
  }
    \toprule
    OD group & Route & Directed node sequence
      & \shortstack{Free-flow\\time (h)} & Familiarity
      & \shortstack{Uses\\$(32,34)$} \\
    \midrule
    $(30,23)$ & $p^0$ & $30\to31\to32\to34\to33\to24\to23$ & $0.395626$ & $0.950000$ & Yes \\
              & $p^1$ & $30\to31\to32\to33\to24\to23$            & $0.397317$ & $0.789060$ & No \\
              & $p^2$ & $30\to31\to23$                              & $0.499182$ & $0.572415$ & No \\
    \addlinespace[0.3em]
    $(30,24)$ & $p^0$ & $30\to31\to32\to34\to33\to24$            & $0.359837$ & $0.950000$ & Yes \\
              & $p^1$ & $30\to31\to32\to33\to24$                  & $0.361528$ & $0.788966$ & No \\
              & $p^2$ & $30\to60\to32\to34\to33\to24$            & $0.490145$ & $0.550331$ & Yes \\
    \addlinespace[0.3em]
    $(31,23)$ & $p^0$ & $31\to32\to34\to33\to24\to23$            & $0.363267$ & $0.950000$ & Yes \\
              & $p^1$ & $31\to32\to33\to24\to23$                  & $0.364958$ & $0.788976$ & No \\
              & $p^2$ & $31\to23$                                    & $0.466823$ & $0.567285$ & No \\
    \addlinespace[0.3em]
    $(31,24)$ & $p^0$ & $31\to32\to34\to33\to24$                  & $0.327478$ & $0.950000$ & Yes \\
              & $p^1$ & $31\to32\to33\to24$                        & $0.329169$ & $0.788864$ & No \\
              & $p^2$ & $31\to32\to34\to33\to27\to26\to24$      & $0.459698$ & $0.481174$ & Yes \\
    \bottomrule
  \end{tabularx}
  \endgroup
\end{table}

All defaults and the second alternatives for $(30,24)$ and $(31,24)$ use arc $(32,34)$, where incidents multiply free-flow time by $3.5$ and capacity by $0.2$. Incidents affect only this arc, and all routes remain open. All arcs use BPR parameters $(\alpha_{\mathrm{BPR}},\beta_{\mathrm{BPR}})=(0.15,4)$ and capacities in travelers per $0.25$-hour period.

\subsubsection{Demand generation, traffic propagation, and cost accounting.}

The EMA trip table gives the same demand of $658.05$ for each of the four selected OD pairs. To stress the shared bottleneck, we double each demand to $1316.10$ and place $80\%$ in the initial queue, giving $q_{o,0}=1053$ after rounding. The network starts empty. For $t=0,\ldots,11$, departures follow $D^{\mathrm{act}}_{o,t}\mid q_{o,t} \sim\operatorname{Binomial}(q_{o,t},\rho_t^{\mathrm{dep}})$, where
\[
    (\rho_t^{\mathrm{dep}})_{t=0}^{11}
    =(0.20,0.25,0.30,0.35,0.40,0.45,0.50,0.55,0.62,0.70,0.80,0.90).
\]
The remaining $20\%$ arrives as independent Poisson counts with means
$\mathbb E[D^{\mathrm{new}}_{o,t+1}]=263.22w_t$, where
\[
    (w_t)_{t=0}^{11}
    =(0.14,0.14,0.13,0.12,0.10,0.09,0.08,0.07,0.05,0.04,0.03,0.01).
\]
New travelers join the queues after period $t$'s departures and traffic loading. Travelers who join after period $11$ remain in the terminal queues. Because the departure probabilities and arrival means depend only on $t$, the demand kernel $K_t^D(\cdot\mid\mathbf q_t)$ depends only on the period and the current queues.

Traffic propagation follows Section~\ref{subsec:dp-traffic-cost} with $\Delta t=0.25$ hours. A packet with an arc duration of one period can enter its next arc in the following period without additional waiting. The simulator uses unconstrained link storage and omits spillback.

The stage cost~\eqref{eq:dp-stage-cost} charges each entering traveler the full unrounded arc time, including traversals that end after $T=12$, and the terminal penalty is the one stated in Section~\ref{subsec:exp-setup}. We report the sum of arc-entry costs and terminal penalties before reward transformations.

\subsubsection{Road observations, compliance filtering, and performance updating.}
\label{app:ema-beliefs}

The road state is normal or incident. The incident probability is $0.30$ in period $0$, $0.18$ after a normal period, and $0.75$ after an incident period. Policies that update road beliefs observe a period-end alarm with sensitivity $0.85$ and false-positive probability $0.10$, along with error-free times on all links. They apply~\eqref{eq:road-filter-update} to both observations and the context in~\eqref{eq:road-filter-context}, then apply~\eqref{eq:road-filter-predict} to predict the next road state. The alarm can inform the joint filter when link times leave the state uncertain. In E2 and E3, given the realized entry flows, the observed time on arc $(32,34)$ alone identifies the completed period's road state. The next-decision incident probability is therefore $0.18$ or $0.75$. The current-road oracle observes $\Theta_t$ before acting.

We compute pre-decision features from $b_t^\Theta$ and inherited entries $h_{ij,t}(\mathbf X_t)$ before the period's new departures. The belief-weighted arc parameters are $\widehat t^0_{ij,t} = \sum_{\vartheta\in\mathcal S^\Theta} b_t^\Theta(\vartheta)t^0_{ij,t}(\vartheta)$ and $\widehat c_{ij,t}= \sum_{\vartheta\in\mathcal S^\Theta} b_t^\Theta(\vartheta)c_{ij,t}(\vartheta)$, in hours and travelers per $0.25$-hour period, respectively.

The unrounded BPR route predictions are 
\begin{equation}
\label{eq:app-ema-predicted-route-time}
\widehat{\operatorname{TT}}^{\mathrm{pred}}_{p,t}
=
\sum_{(i,j)\in p}
\widehat t^0_{ij,t}
\left[
1+\alpha_{\mathrm{BPR}}
\left(
\frac{h_{ij,t}(\mathbf X_t)}{\widehat c_{ij,t}}
        \right)^{\beta_{\mathrm{BPR}}}
\right].
\end{equation}
Recommended and default routes use the same arc-level predictions. We set $s_{\mathrm{TT}}=0.25$ hours in~\eqref{eq:compliance-travel-time-difference} and $C_{\mathrm{Cong}}=5$ in~\eqref{eq:compliance-congestion-feature}.

We leave $\Delta\mathrm{TT}_{op,t}$ unclipped and fix both features during likelihood evaluation. Because BPR is nonlinear, predictions from belief-weighted parameters generally differ from posterior averages of state-specific route times at the same inherited flows.

Each trajectory with uncertain compliance starts with $M=3000$ particles that we draw independently from the diagonal Gaussian prior and weight equally at $1/M$. E2 uses the biased prior in Section~\ref{subsec:exp-e2}. E3 adds the accurate and diffuse priors in Section~\ref{subsec:exp-e3}. In E3, common scenario seeds give frozen and full-adaptation policies identical initial particles within each prior and training replicate. Although the prior coordinates are independent, the weighted particles represent the joint posterior, which need not factor across coordinates.

When a policy updates compliance beliefs, we evaluate $L_t^{C,M}$ (Section~\ref{subsec:particle-approximation}) using~\eqref{eq:compliance-likelihood}, the realized integer counts, and the features that we compute from $S_t^{(M)}$ before acting. These features are the same for every particle. Only informative cells with positive counts contribute non-constant factors. We apply~\eqref{eq:particle-weight-update} in log space and clip Bernoulli probabilities to $[10^{-12},1-10^{-12}]$ before taking logarithms. When $\operatorname{ESS}_{t+1}<1500$, we resample systematically and reset the weights to $1/M$. Without informative observations, the compliance belief stays unchanged. Policies with oracle compliance information set all $3000$ particles to the true vector and use a resampling threshold of zero.

The performance signal~\eqref{eq:performance-signal} weights informative routes by their realized recommendation counts and compares them with their default routes on the same period-end snapshot. Frozen policies fix compliance particles, road beliefs, and performance states at their initial values while traffic and queues evolve. We hold actor parameters fixed during validation and evaluation.

\subsubsection{Policy input and PPO implementation.}
\label{app:ema-ppo-implementation}

The policy input $\widetilde S_t=\psi(S_t^{(M)})$ has $315$ coordinates, which we normalize before passing them to the actor and critic (Section~\ref{subsec:state-encoding}). It contains the decision period index ($1$), waiting traveler counts by group ($4$), the inherited entry count on each directed link ($258$), the total number of unfinished travelers ($1$), the pre-decision incident probability ($1$), recommendation performance states by group ($4$), posterior means and marginal posterior variances of the compliance coefficients ($5$ each), and route-level predictive compliance probabilities, predicted travel times, and availability indicators ($12$ each). The $259$ traffic coordinates consist of the $258$ inherited entry counts $h_{ij,t}(\mathbf X_t)$ and $U(\mathbf X_t)$. They omit packet routes and residual durations. Posterior means and marginal variances summarize the weighted particles without recording posterior correlations. We use a fixed order over all $258$ links. Because the network starts empty, only the $13$ links in the route catalog can have nonzero inherited entry counts.

All twelve route-availability indicators equal one before normalization, and the policy input has no separate active-group indicators. The actor produces shares $\mathbf r_t^{\mathrm{actor}}$ for all four groups. All four Dirichlet factors in~\eqref{eq:joint-dirichlet-policy} enter the likelihood ratio~\eqref{eq:ppo-ratio} and the entropy term in~\eqref{eq:ppo-full-objective}. The action $\mathbf r_t$ retains the actor's shares for groups with positive queues, as in Section~\ref{subsec:ppo-actor}. The simulator uses those shares when travelers depart.

The actor and critic share two fully connected $128$-unit hidden layers with hyperbolic-tangent activations and use separate output layers. The actor's twelve outputs form four three-route groups, with concentrations from~\eqref{eq:dirichlet-concentration} and $\kappa_{\min}=0.05$. The critic outputs one value. Orthogonal initialization uses gains $\sqrt{2}$, $0.01$, and $1$ for the hidden layers, actor output layer, and critic output layer, respectively, with zero biases. Table~\ref{tab:app-ema-ppo-settings} lists the remaining settings.

\begin{table}[htbp!]
  \centering
  \caption{PPO configuration for each E2--E3 actor.}
  \label{tab:app-ema-ppo-settings}
  \begingroup
  \small
  \setlength{\tabcolsep}{5pt}
  \renewcommand{\arraystretch}{1.10}
  \begin{tabularx}{\linewidth}{
    @{}>{\raggedright\arraybackslash}Xr
       >{\raggedright\arraybackslash}Xr@{}
  }
    \toprule
    Setting & Value & Setting & Value \\
    \midrule
    Training iterations & 800 & Transitions per iteration & 960 \\
    Complete trajectories per iteration & 80 & Update epochs & 4 \\
    Minibatch size & 320 & Gradient-norm limit & 0.50 \\
    Adam learning rate & $3\times10^{-4}$ & Adam numerical constant & $10^{-5}$ \\
    Discount $\gamma$ & 1.00 & GAE parameter $\lambda$ & 0.95 \\
    PPO clip radius & 0.20 & Entropy coefficient & 0.01 \\
    Value-loss coefficient & 0.50 & Minimum concentration & 0.05 \\
    Reward scale & 0.01 & Reward clip & 10.0 \\
    Observation clip & 10.0 & Running-statistic pseudocount & $10^{-4}$ \\
    Action-share floor & 0 & Device & CPU \\
    \bottomrule
  \end{tabularx}
  \endgroup
\end{table}

During training, the actor samples Dirichlet shares. Validation and evaluation use the mean shares~\eqref{eq:dirichlet-mean-action}, while demand, assignments, compliance, and road conditions remain stochastic. For numerical stability, the implementation clips shares below $10^{-8}$ and renormalizes each group's vector before computing log densities.

We standardize observations with running means and variances and clip them to $[-10,10]$. We scale rewards by $0.01$, divide them by their running standard deviation without centering, and clip them to $[-10,10]$. Residuals and advantages~\eqref{eq:td-residual}--\eqref{eq:gae} use transformed rewards and the stored critic values, and we standardize advantages once per batch. Only training updates the normalization statistics, which we save with each checkpoint and reuse in validation and evaluation. Each actor uses $64{,}000$ training trajectories ($768{,}000$ transitions) and $9600$ minibatch updates.

\subsubsection{Out-of-sample evaluation and penalties for travelers still en route.}
\label{app:ema-evaluation}

E2 evaluates three independently trained PPO actors under the biased prior and nine benchmark policies on $1000$ common out-of-sample scenarios (Table~\ref{tab:exp-e2-ema}). The predicted-time softmax assigns shares within each group proportional to $\exp(-\widehat{\operatorname{TT}}^{\mathrm{pred}}_{p,t}/\lambda_{\mathrm{sm}})$, with $\lambda_{\mathrm{sm}}=0.40$ hours. The route times in~\eqref{eq:app-ema-predicted-route-time} use inherited arc entries and the current road belief before new departures.

The route-aware myopic benchmarks start coordinate descent from the default-only allocation and use up to three sweeps with $0.2$ share increments. For each candidate allocation, they multiply each group's expected departures $q_{o,t}\rho_t^{\mathrm{dep}}$ by its route shares and split the resulting counts between each recommended route and its default using the benchmark's compliance probabilities. Let $m_{op,t}$ denote the resulting expected number of travelers who depart on route $p$ in group $o$. The proxy adds these travelers to every arc of their route and includes inherited entries, giving arc flows
\[
    \widehat z_{ij,t}
    =h_{ij,t}(\mathbf X_t)
      +\sum_{o\in\mathcal O}\sum_{p\in\mathcal P_{o,t}(S_t):(i,j)\in p} m_{op,t}.
\]
Let $\widehat t_{ij,t}(z)$ denote the unrounded BPR time at flow $z$ with the belief-weighted parameters of Section~\ref{app:ema-beliefs}. The road oracle uses the parameters of the observed current state. The full-route proxy is
\begin{equation}
\label{eq:app-ema-myopic-proxy}
    \sum_{(i,j)\in\mathcal E}\widehat z_{ij,t}\widehat t_{ij,t}(\widehat z_{ij,t})
    + 5\sum_{o\in\mathcal O}\sum_{p\in\mathcal P_{o,t}(S_t)}
       m_{op,t}\left(\sum_{(i,j)\in p}
         \left\lceil\widehat t_{ij,t}(\widehat z_{ij,t})/\Delta t\right\rceil-1\right).
\end{equation}
The first term charges all proxy flows at their unrounded arc times. The second charges $5$ traveler-hours per expected departing traveler per predicted traversal period beyond the current period. Section~\ref{subsec:exp-setup} specifies the information and allocation updates for each benchmark.

\begin{figure}[htbp!]
  \centering
  \includegraphics[width=0.85\linewidth]{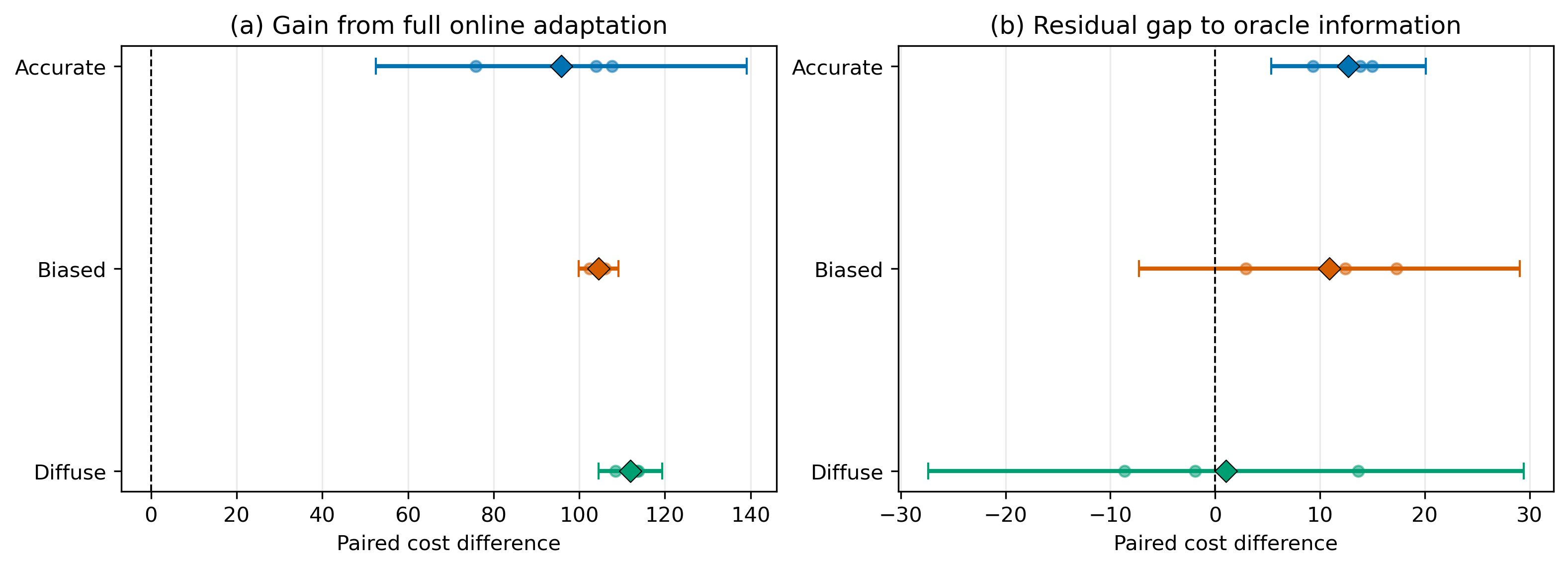}
  \caption{E3 paired mean differences across training replicates. Panel~(a) reports \emph{frozen} cost minus \emph{full adaptation} cost. Panel~(b) reports \emph{full adaptation} cost minus \emph{oracle information} cost. Translucent circles show the three training-replicate differences, diamonds show their means, and whiskers show $95\%$ Student-$t$ intervals across replicates. 
  Positive values in Panel~(a) favor full adaptation. A zero difference in Panel~(b) denotes equal mean costs for that comparison with the separately trained oracle PPO policy.}
  \label{fig:exp-e3-paired-effects}
\end{figure}

For each prior and information setting, E3 uses three independently trained actors. The biased-prior full-adaptation setting uses the same actors as E2. We validate the mean policy on $40$ held-out scenarios after the first iteration and every $25$ iterations. We select the minimum-cost checkpoints and evaluate them with their saved normalization statistics. Within each prior and training replicate, the three information settings share $300$ out-of-sample scenarios. We use a distinct set for each prior and training replicate. Figure~\ref{fig:exp-e3-paired-effects} shows the replicate-level paired differences that underlie Table~\ref{tab:exp-e3-adaptation}.

Section~\ref{subsec:exp-setup} defines the cost metrics, the subtraction order of each paired comparison, and the Student-$t$ intervals across training replicates. Scenario-level intervals measure evaluation uncertainty conditional on the selected actors.

\subsection{Sensitivity to the Penalty for Travelers Still En Route}
\label{app:terminal-sensitivity}

Because terminal penalties drive most of the cost differences in E2 and E3, we examine whether the penalty for travelers still en route drives the policy comparisons. We recompute saved out-of-sample costs as
\begin{equation}
\label{eq:app-terminal-revaluation}
    \operatorname{Cost}(\eta_X)
    =
    \sum_{t=0}^{11}g_t(S_t,\mathbf r_t,\Xi_{t+1})
    +20\sum_{o\in\mathcal O}q_{o,12}
    +\eta_X U(\mathbf X_{12}),
    \quad \eta_X\in\{0,2.5,5,10\}.
\end{equation}
The penalty $\eta_X$ is the cost in traveler-hours for each traveler still en route after three hours. The primary experiments use $\eta_X=5$. We keep the saved recommendations, responses, road conditions, traffic movements, and terminal counts and recompute the paired mean differences and intervals without retraining or making new routing decisions. At $\eta_X=0$, the paired differences reflect arc-entry costs alone because the policies in each pair have identical terminal queues and their waiting penalties cancel.

Table~\ref{tab:exp-terminal-sensitivity} reports the resulting paired mean cost differences. In E2, PPO retains a lower mean arc-entry travel cost than the oracle-compliance myopic benchmark, with an interval that excludes zero. In E3, the point estimates favor full adaptation over the frozen policy, although the accurate- and diffuse-prior intervals include zero. At half the primary penalty ($\eta_X=2.5$), all three E3 intervals favor full adaptation. Every point estimate favors PPO or full adaptation at all four penalty values, and the advantages grow as $\eta_X$ increases. A planner who places more weight on clearing the network therefore gains more from dynamic, adaptive recommendations, and a planner who cares only about travel time still benefits from full-horizon planning.

\begin{table}[htbp!]
  \centering
  \caption{Cost sensitivity to the penalty per traveler still en route after three hours. We recalculate costs for the same out-of-sample simulation runs without changing policies or traffic outcomes.  
  The columns give paired mean cost differences in hours with $95\%$ Student-$t$ intervals across training replicates. Column ``E2'' subtracts PPO cost from \emph{oracle-compliance myopic cost}. Columns ``E3 accurate'', ``E3 biased'', and ``E3 diffuse'' subtract \emph{full adaptation} cost from \emph{frozen} cost. Positive values favor PPO or full adaptation.}
  \label{tab:exp-terminal-sensitivity}
  \scriptsize
  \setlength{\tabcolsep}{4pt}
  \begin{tabular}{rrrrr}
    \toprule
    $\eta_X$ & E2 & E3 accurate & E3 biased & E3 diffuse \\
    \midrule
    $0.0$  & $55.4\ [53.3,57.5]$   & $11.2\ [-7.3,29.7]$  & $8.4\ [7.6,9.2]$     & $5.8\ [-3.1,14.7]$ \\
    $2.5$  & $121.2\ [117.0,125.4]$& $53.5\ [22.9,84.1]$  & $56.5\ [53.7,59.2]$  & $58.8\ [51.9,65.8]$ \\
    $5.0$  & $187.0\ [180.6,193.3]$& $95.8\ [52.6,139.1]$ & $104.5\ [99.9,109.2]$& $111.9\ [104.5,119.4]$ \\
    $10.0$ & $318.5\ [307.9,329.1]$& $180.5\ [111.3,249.6]$& $200.6\ [192.1,209.2]$& $218.0\ [204.2,231.8]$ \\
    \bottomrule
  \end{tabular}
\end{table}

\end{appendices}

\end{document}